\documentclass{article}
   \usepackage{amssymb} 
 \usepackage[totalwidth=410pt, totalheight=620pt]{geometry}
\newcommand{\g}{\mbox{$\bf g$}}

\newcommand{\h}{\mbox{\textbf{h}}}
\newcommand{\n}{\mbox{\textbf{n}}}

\newcommand{\al}{\alpha}
\newcommand{\eps}{\epsilon}
\newcommand{\la}{\lambda}
\newcommand{\La}{\Lambda}

\newcommand{\mb}{\mbox}

\newcommand{\Mklz}[2]{\left\{\left.\;#1\;\right|\; #2\;\right\}}
\newcommand{\F}{\mathbb{F}}

\newcommand{\N}{\mathbb{N}}
\newcommand{\Nn}{\mathbb{N}_0}
\newcommand{\Q}{\mathbb{Q}}

\newcommand{\R}{\mathbb{R}}

\newcommand{\Z}{\mathbb{Z}}
\newcommand{\FK}[1]{\mbox{$\F\,[#1]$}}

\newcommand{\Fi}[1]{{\mathcal F}_{#1}}
\newcommand{\FA}{\mbox{${\mathcal F}_A$}}

\newcommand{\FC}{\mbox{${\mathcal F}_C$}}
\newcommand{\FD}{\mbox{${\mathcal F}_D$}}

\newcommand{\Fa}[1]{ Fa(#1) }
\newcommand{\Proof}{\mbox{\bf Proof: }}
\newcommand{\qed}{\mb{}\hfill\mb{$\square$}\\}
\newcommand{\End}{\mb{}\hfill\mb{$\square$}\\}
\newcommand{\Spm}{\mbox{Specm\,}}

\newcommand{\ti}{\tilde}
\newcommand{\res}[1]{\!\mid_{#1}}

\newcommand{\ro}{{\,\bf \triangleright\,}}
\newcommand{\lo}{{\,\bf \triangleleft\,}}
\begin{document}
\newtheorem{Theorem}{Theorem}[section]
\newtheorem{Definition}[Theorem]{Definition}
\newtheorem{Proposition}[Theorem]{Proposition}
\newtheorem{Corollary}[Theorem]{Corollary}
\newtheorem{Conjecture}[Theorem]{Conjecture}
\newtheorem{Remark}[Theorem]{Remark}
\newtheorem{Remarks}[Theorem]{Remarks}
\newtheorem{Lemma}[Theorem]{Lemma}
\newtheorem{Example}[Theorem]{Example}
\newtheorem{Examples}[Theorem]{Examples}
%
%
\title{Integrating infinite-dimensional Lie algebras by a Tannaka
  reconstruction (Part I)}
\author{Claus Mokler\thanks{Supported by the Deutsche Forschungsgemeinschaft.}\\\\ Universit\"at Wuppertal, Fachbereich C - Mathematik\\  Gau\ss stra\ss e 20\\ D-42097 Wuppertal, Germany\vspace*{1ex}\\ 
          mokler@math.uni-wuppertal.de}
\date{}
\maketitle
\begin{abstract}\noindent 
Let $\g$ be a Lie algebra over a field $\F$ of characteristic zero, let $\mathcal C$ be a certain tensor category of representations of $\g$, and
${\mathcal C}^{du}$ a certain category of duals. By a Tannaka reconstruction we associate to $\mathcal C$ and ${\mathcal C}^{du}$ a monoid
$M$ with a coordinate ring of matrix coefficients $\FK{M}$ (which has in general no natural coalgebra structure), as well as a Lie algebra $Lie(M)$. 
The monoid $M$ and the Lie algebra $Lie(M)$ both act on the objects of $\mathcal C$.  
We say that $\mathcal C$ and ${\mathcal C}^{du}$ are good for integrating $\g$, if the Lie algebra $\g$ is in a natural way a Lie subalgebra of 
$Lie(M)$. In this situation we treat:
The adjoint action of the unit group of $M$ on $Lie(M)$, the relation between $\g$ and $M$-invariant subspaces, 
the embedding of the coordinate ring $\FK{M}$ into the dual of the
universal enveloping algebra $U(\g)$ of $\g$, a Peter-and-Weyl-type theorem
for $\FK{M}$ if $\mathcal C$ is a semisimple category, the Jordan-Chevalley decompositions for
certain elements of $M$ and $Lie(M)$, an embedding theorem related to
subalgebras of $\g$ which act locally finite, prounipotent subgroups and generalized toric submonoids of $M$.
We show that $\mathcal C$ and ${\mathcal C}^{du}$ are good for integrating $\g$ if the Lie
algebra $\g$ is generated by integrable locally finite elements. 

We interprete the monoid $M$ algebraic geometrically as an irreducible  weak algebraic monoid with
Lie algebra $Lie(M)$. The monoid $M$ acts by morphisms of varieties on every
object $V$ of $\mathcal C$. The action of the Lie algebra $Lie(M)$ on $V$ is the differentiated action. 
\end{abstract}
\vspace*{1ex} 
\section*{Introduction}
Finite dimensional Lie algebras arise as an important aid in the investigation of finite dimensional Lie groups and algebraic groups, 
encoding their local structure. 
In the infinite-dimensional situation one often finds natural 
infinite-dimensional Lie algebras, which also have a rich representation theory. But the Lie algebras are obtained without 
the help of any groups. One faces the problem to ``integrate infinite dimensional Lie algebras'', i.e., to construct groups 
associated in an appropriate way to  these Lie algebras.   
Maybe the most important examples are Kac-Moody algebras and their associated Kac-Moody groups. This example also shows that one finds, in 
difference to the classical situation, many different but closely related groups, like minimal Kac-Moody groups, formal Kac-Moody groups, etc.
Furthermore, there are related algebraic geometric constructions, the most important the flag varieties. 
There exists no general theory of integrating infinite-dimensional Lie algebras in an algebraic geometric context. The considerations given 
in this paper and in \cite{M4} indicate, that it may be possible to base a theory on the Tannaka reconstruction, 
and on an algebraic geometry, which is to some extend motivated by the algebraic geometry of affine ind-varieties of 
I. R. \v{S}afarevi\v{c}, \cite{Sa}.  But much work remains to be done.\vspace*{1ex}

As a particular case of Tannaka-Krein theory, a linear algebraic group $G$ over a field $\F$, and its coordinate ring $\FK{G}$ can be 
reconstructed from the $\F$-linear abelian symmetric rigid tensor category of
rational $G$-modules. Compare for example the article \cite{DeMi}, the survey \cite{JS}, or Section
2.5 of the book \cite{Sp}. Now every 
rational $G$-module differentiates to a module of its Lie algebra $Lie(G)$, and every $G$-equivariant homomorphism between two $G$-modules 
is also $Lie(G)$-equivariant. For an irreducible linear algebraic group $G$ over a field $\F$ of 
characteristic zero the category of rational $G$-modules is isomorphic in this way to the corresponding
category of $Lie(G)$-modules. Therefore, we may start the reconstruction also
with the corresponding category of $Lie(G)$-modules. Note that it contains with every
$Lie(G)$-module $V$ also the dual $Lie(G)$-module $V^*$. 
A similar reconstruction holds for an irreducible linear algebraic monoid $M$ over a field $\F$ of characteristic zero and its coordinate ring 
$\FK{M}$, where we have to start with the appropriate category of representations of the Lie
algebra $Lie(M^\times)$ of the unit group $M^\times$ of $M$, which is an
irreducible linear algebraic group. If the $Lie(M^\times)$-module $V$ is
contained in this category, then its dual $V^*$ does not need to be contained
in the category.   

If we want to integrate infinite-dimensional Lie algebras over fields of characteristic zero by a Tannaka reconstruction we have to generalize 
the reconstruction of irreducible linear algebraic monoids, which includes the reconstruction of irreducible linear algebraic groups. 
The reason for this is the role of the duals. For an infinite-dimensional representation of a Lie algebra $\g$ there may be several natural 
possibilities to choose a dual representation. For example for an infinite-dimensional integrable highest weight representations of a Kac-Moody 
algebra the integrable or restricted dual, as well as the full dual are important possibilities.
In general, the duals don't fit in a tensor category, since the tensor product of the duals of infinite-dimensional 
representations does not need to be isomorphic in the natural way to the dual of the tensor product of these representations. 
In general there is no rigidity any more, a double dual of an infinite-dimensional representation may be much bigger than the representation 
itself. 
Linear algebraic monoids are interesting objects with a rich theory. For an
excellent survey compare \cite{Re}. Similarly, the monoids obtained by integrating infinite-dimensional Lie
algebras are interesting. For an example compare \cite{M1}, \cite{M2}, and \cite{M3}. \vspace*{1ex}

In Section \ref{TK} we start with a Lie algebra $\g$ (of arbitrary dimension) over a field $\F$ of characteristic zero,  
a certain category $\mathcal C$ of representations of $\g$, and a certain category of duals ${\mathcal C}^{du}$. 
By a Tannaka reconstruction we associate to $\mathcal C$ and ${\mathcal C}^{du}$ the Tannaka monoid $M$ with coordinate ring of matrix 
coefficients $\FK{M}$, as well as a Lie algebra $Lie(M)$. The monoid $M$ and the Lie algebra $Lie(M)$ both act on the objects of $\mathcal C$.  
We call $\mathcal C$ and ${\mathcal C}^{du}$ good for integrating $\g$, if the Lie algebra $\g$ identifies in a natural way with a Lie
subalgebra of $Lie(M)$.

The important thing here is to work with two categories. The coordinate rings which we obtain in this way have a weaker algebraic structure as it 
is usual from other Tannaka reconstruction situations. In general, the multiplication of the Tannaka monoid $M$ does not induce a
comultiplication of its coordinate ring $\FK{M}$. Only the left and right multiplications with elements of $M$ induce comorphisms.

We investigate the Tannaka monoid $M$, its coordinate ring of matrix coefficients $\FK{M}$, and the Lie algebra $Lie(M)$ in the case where 
$\mathcal C$ and ${\mathcal C}^{du}$ are good for integrating $\g$. 
We determine the kernel of the adjoint action of the unit group $M^\times$ on the Lie algebra $Lie(M)$. We investigate the relation between 
the $\g$-invariant subspaces and the $M$-invariant subspaces of the $\g$-modules contained in $\mathcal C$. We show that the coordinate ring 
$\FK{M}$ embeds as an algebra into the dual of the universal enveloping algebra $U(\g)$ of $\g$. In particular, this implies 
that $M$ is irreducible. If the category $\mathcal C$ is semisimple, (i.e., every $\g$-module contained in $\mathcal C$ is a possibly 
infinite sum of irreducible $\g$-modules contained in $\mathcal C$), we show a Peter-and-Weyl-type theorem for the coordinate ring
$\FK{M}$. We give an additive Jordan-Chevalley decomposition for certain elements of $Lie(M)$, and a multiplicative Jordan-Chevalley decomposition 
for certain elements of certain groups associated to the idempotents of $M$. 
We show an embedding theorem related to subalgebras of $\g$ which act locally finite on every $\g$-module contained in $\mathcal C$. 
In particular, if a subalgebra of $\g$ acts locally nilpotent on the modules contained in $\mathcal C$, and if a compatibility condition with the
duals is satisfied, we get a prounipotent subgroup of $M$. If the $\g$-modules contained in $\mathcal C$ have weight space
decompositions with respect to some abelian subalgebra of $\g$, if the set of weights are contained in a lattice, 
and if a compatibility condition with the duals is satisfied, we get a generalized toric submonoid of $M$. 

In general it is a difficult problem to decide if some categories $\mathcal C$ and ${\mathcal C}^{du}$ are good for integrating $\g$.
We obtain the following criterion, which is to expect to hold if the Tannaka reconstruction is really a reasonable way 
to integrate infinite-dimensional Lie algebras:
If the Lie algebra $\g$ is generated by integrable locally finite elements, then $\mathcal C$ and $\mathcal C^{du}$ are good for integrating $\g$. 
Furthermore, there exists a dense submonoid of $M$, build in a similar manner as the (minimal) Kac-Moody group of \cite{KP1}, 
or as the associated groups of linear algebraic integrable Lie algebras of \cite{K1}.\vspace*{1ex} 

In Section \ref{AGeo} we interprete the monoid $M$ algebraic geometrically as an irreducible weak algebraic monoid with Lie algebra $Lie(M)$. 
The monoid $M$ acts by morphisms of varieties on every object $V$ of $\mathcal C$. The action of the Lie algebra $Lie(M)$ on $V$ is 
the differentiated action.\vspace*{1ex}

In the subsequent paper \cite{M4} we introduce regular (or polynomial) functions on a dense subgroup of $M$, as well as linear regular functions on 
$U(\g)$ in the case where $\g$ is generated by one-parameter elements. We investigate and describe various coordinate rings of matrix coefficients 
associated to integrable representations of a Lie algebra. 
This is motivated by some coordinate rings for Kac-Moody groups in \cite{KP2}, \cite{Kas}, by the coordinate rings of the associated
groups of linear algebraic integrable Lie algebras in \cite{K1}, and by the
polynomial functions on free Kac-Moody groups and some results related to the shuffle algebra in \cite{BiPi}.

Even if we know that $\mathcal C$ and ${\mathcal C}^{du}$ are good for integrating an infinite-dimensional Lie algebra, it may be quite hard, 
and also a long way to determine the associated Tannaka monoid and its Lie algebra explicitely.  
An example has been treated in \cite{M1}: The category of integrable $\g$-modules contained in the category $\mathcal O$ of a symmetrizable 
Kac-Moody algebra $\g$ is one possible generalization of the category of finite-dimensional representations of a semisimple Lie algebra. 
It keeps  the complete reducibility theorem, i.e., every integrable $\g$-module contained in $\mathcal O$ is a sum of integrable irreducible 
highest weight modules. In \cite{M1} we determined the Tannaka monoid associated to this category and its category of integrable duals. 
Its Zariski-open dense unit group  is the (minimal) Kac-Moody group. Its Lie algebra identifies with the Kac-Moody algebra.
To treat this example did not require an investigation of the Tannaka reconstruction since the Kac-Moody group, and the restriction of the 
algebra of matrix coefficients onto the Kac-Moody group had been well known by investigations of V. Kac, D. Peterson and by M. Kashiwara.
In \cite{M4} we also determine the Tannaka monoid associated to the category of integrable $\g$-modules contained in the category $\mathcal O$ 
and its category of full duals. Its Zariski-open dense unit group is the formal Kac-Moody group. Its Lie algebra identifies with the formal 
Kac-Moody algebra.\vspace*{1ex}

Often it is also possible to integrate infinite-dimensional Lie algebras by a functional analytical version of the Tannaka reconstruction, 
even if the algebraic geometric integration fails. In a subsequent paper \cite{M6} we treat as an example generalized Kac-Moody algebras.

\tableofcontents
%
%
\section{Preliminaries}\label{Pre}
%
%
%

In this paper we denote by $\N=\Z^+$, $\Q^+$, resp. $\R^+$ the sets of strictly positive numbers of $\Z$, $\Q$, resp. $\R\,$,
and the sets $\N_0=\Z^+_0$, $\Q^+_0$, $\R^+_0$ contain in addition the zero.
In the whole paper $\F$ is a field of characteristic 0, and $\F^\times$ denotes its unit group.

\subsection*{The category of sets with coordinate rings}
A {\it set with coordinate ring} consists of a non-empty
set $A$ equipped with a point separating algebra of functions $\FK{A}$, which
we call {\it coordinate ring}.
We equip $A$ with the {\it Zariski topology}, i.e., the closed sets are given by the
common zero sets of the functions of $\FK{A}$. Note that $A$ is irreducible 
if and only if $\FK{A}$ is an integral domain.

A {\it morphism} from the set with coordinate ring ($A,\FK{A}$) to the set with
coordinate ring ($B,\FK{B}$) consists of a map $\phi:A\to B$, whose comorphism 
$\phi^*:\FK{B}\to\FK{A}$ exists. In particular a morphism is Zariski continuous.

If ($B,\FK{B}$) is a set with coordinate ring, and $A$ is a non-empty subset of
$B$, then we get a coordinate ring on $A$ by restricting the 
functions of $\FK{B}$ to $A$.

If ($A,\FK{A}$) is a set with coordinate ring and $f\in\FK{A}\setminus\{0\}$,
the principal open set $D_A(f):=\Mklz{a\in A}{f(a)\neq 0}$ is equipped with a
coordinate ring by identifying the localization $\FK{A}_f$ in the obvious way
with an algebra of functions on $D_A(f)$. If $A$ is irreducible, then also
$D_A(f)$ is irreducible.

If ($A,\FK{A}$) and ($B,\FK{B}$) are sets with coordinate rings, then the
product $A\times B$ is equipped with a coordinate ring by identifying the tensor product 
$\FK{A}\otimes\FK{B}$ in the obvious way with an algebra of functions
on $A \times B$. If $A$ and $B$ are irreducible, then also $A\times B$ is irreducible.

Let ($A,\FK{A}$) and ($B,\FK{B}$) sets with coordinate rings. A morphism
$\phi: A\to B$ is called an embedding if the restriction $\phi:A\to\phi(A)$ is
an isomorphism of sets with coordinate rings. An embedding $\phi: A\to B$ is
called closed if its image $\phi(A)$ is closed in $B$.
\vspace*{1ex}

For a set with coordinate ring ($A,\FK{A}$) and $a\in A$ we denote by
$Der_a(\FK{A})$ the $\F$-linear space of derivations in $a$.
We denote by $Der(\FK{A})$ the Lie algebra of derivations of the algebra $\FK{A}$.
We denote by $\Spm \FK{A}$ the $\F$-valued points of the algebra $\FK{A}$, i.e., the
homomorphisms of algebras from $\FK{A}$ to $\F$.

\subsection*{The dual of the universal enveloping algebra of a Lie algebra}
Let $\g$ be a Lie algebra over $\F$. Equip the universal enveloping algebra $U(\g)$ of $\g$ with its canonical
Hopf algebra structure: The comultiplication $\Delta:\,U(\g)\to U(\g)\otimes
U(\g)$ is the algebra homomorphism induced by $\Delta(x)=1\otimes x+ x\otimes
1$, $x\in\g$. The counit 
$\eps:\,U(\g)\to\F$ is the algebra homomorphism induced by $\eps(x)=0$, $x\in
\g$. The antipode $S:\,U(\g)\to U(\g)$ is the algebra anti-automorphism induced by
$S(x)=-x$, $x\in\g$.

The full dual $U(\g)^*$ is a commutative algebra with unit $\eps$, the multiplication defined by 
\begin{eqnarray*}
   (h\cdot h') (x) := \left(h\otimes h'\right)(\Delta(x))  \quad\mb{ where }\quad h,\,h'\in U(\g)^*,\;x\in U(\g).
\end{eqnarray*} 

\begin{Theorem}\label{zd} Let $\g$ be a Lie algebra over a field $\F$ of characteristic
  zero. Then $U(\g)^*$ has no zero divisors.
\end{Theorem}
\Proof For a finite dimensional
Lie algebra $\g$ the proof can be found in \cite{Ho}, Section XVI.3. We only
give the modifications which are necessary for infinite-dimensional $\g$.

Let $a_i$, $i\in I$, be a base of $\g$. Let $I$ be totally ordered. 
Let $\omega$ be the set of functions $e:I\to\Nn$ such that $e_i:=e(i)\neq 0$ only for
finitely many $i\in I$. We can add elements of $\omega$ by adding the values of the
functions. We denote the zero element by $0$, i.e., the function which maps every element
of $I$ to zero. For $e\in\omega$ define the degree $deg(e):=\sum_{i\in I}e_i$, define the support
$supp(e):=\Mklz{i\in I}{e_i\neq 0}$.

For $a\in U(\g)$ set $a^0:=1$. We get a base of $U(\g)$ by the monomials
\begin{eqnarray*}
 b_e:=\prod_{i\in I}\frac{a_i^{e_i}}{e_i!}\quad\mb{ where }\quad e\in\omega.
\end{eqnarray*}
Here the order of the product is given by the total order of $I$.

For $e\in\omega$ define $h_e\in U(\g)^*$ by $h_e(b_{\ti{e}})=\delta_{e\ti{e}}$, $\ti{e}\in \omega$.
Every element $h$ of $U(\g)^*$ can be identified with an infinite sum
$h=\sum_{e\in\omega}c_e h_e$, $c_e\in\F$. Two elements $h=\sum_{e\in\omega}c_e h_e$ and $\ti{h}=\sum_{e\in\omega}\ti{c}_e h_e$
of $U(\g)^*$ are multiplied by
\begin{eqnarray*}
    h\cdot\ti{h}=\sum_{f\in\omega}\left(\sum_{e,\ti{e}\in\omega,\;e+\ti{e}=f}c_e \ti{c}_{\ti{e}}\right)h_f.
\end{eqnarray*}
(To see this evaluate $h\cdot\ti{h}$ on $b_f$, $f\in\omega$. Prove first by a
direct computation 
\begin{eqnarray*}
  \Delta (b_f)= \sum_{e,\ti{e}\in\omega,\;e+\ti{e}=f} b_e\otimes b_{\ti{e}}\quad \mb{ for all }\quad f\in\omega.)
\end{eqnarray*}

Let $J\subseteq I$ be a finite subset of $I$. For $J\neq\emptyset$ denote by $\F[[x_i\mid i\in J]]$ the algebra of formal power series in the
variables $x_i$, $i\in J$. Set $\F[[x_i\mid i\in \emptyset]]:=\F$.
Assigning to $x_i$ the function $h_i$, $i\in J$, induces an embedding
$\Phi:\,\F[[x_i\mid i\in J]]\to U(\g)^*$ of $\F$-algebras with unit. 

Now suppose that $h,\,\ti{h}\in U(\g)^*$ are both non-zero and
$h\cdot\ti{h}=0$. Write $h$, $\ti{h}$ in the form $h=\sum_{e\in\omega}c_e h_e$, $\ti{h}=\sum_{e\in\omega}\ti{c}_e h_e$.  
Since $h\neq 0$ and $\ti{h}\neq 0$ there exist
coefficients $c_e\neq 0$ and $\ti{c}_{\ti{e}}\neq 0$. Let $J:=supp(e)\cup
supp(\ti{e})$. Write $h$ as the sum
\begin{eqnarray*}
 h=h_J + h_{rest} \quad\mb{ with }\quad
   h_J:=\sum_{f\in\omega,\;supp(f)\subseteq J} c_f h_f\quad \mb{ and }\quad
   h_{rest}:=\sum_{f\in\omega,\;supp(f)\nsubseteq J} c_f h_f.
\end{eqnarray*}
Similarly, $\ti{h}=\ti{h}_J+\ti{h}_{rest}$. We get 
\begin{eqnarray*}
   0=h\cdot\ti{h}= h_J\cdot\ti{h}_J+\left(h_{rest}\cdot\ti{h}_J+ h_J\cdot\ti{h}_{rest}+ h_{rest}\cdot\ti{h}_{rest}\right).
\end{eqnarray*}
It follows $0=h_J\cdot\ti{h}_J$. Since $h_J$ and $\ti{h}_J$ are in the
image of $\Phi$, and $\F[[x_i\mid i\in J]]$ has no zero divisors (to see this apply
Proposition 5.8 (iii) of \cite{H} recursively $|J|$ times), we get
$h_J=\ti{h}_J=0$. This contradicts $c_e\neq 0$ and $\ti{c}_{\ti{e}}\neq 0$.
\qed
\subsection*{Locally finite endomorphisms of a linear space}
All definitions and results which we state in this subsection can be found in,
or follow easily from Section 7.1 of the book \cite{MoPi}. As before $\F$ denotes a field of
characteristic zero, $\overline{\F}$ denotes the algebraic closure of $\F$.\vspace*{1ex}

Let $V$ be a $\F$-linear space. An endomorphism $x\in End(V)$ is called {\it locally finite}, if each $v\in V$ lies in a $x$-invariant
finite dimensional subspace of $V$.
Particular examples of locally finite endomorphisms are diagonalizable,
semisimple, locally nilpotent, locally unipotent, and locally weak unipotent
endomorphisms, which are defined as follows:
$x\in End(V)$ is called {\it diagonalizable} if there exists a base of $V$ consisting of eigenvectors of $x$.
$x\in End(V)$ is called {\it semisimple} if its extension to $End(V\otimes_{\F}\overline{\F})$ is semisimple.
$x\in End(V)$ is called {\it locally nilpotent} if for every $v\in V$ there exists an positive integer $n\in\N$ such that $x^n v=0$.
$x\in Aut(V)$ is called {\it locally unipotent} if $x-id_V$ is locally nilpotent.
$x\in End(V)$ is called {\it weak locally unipotent} if $V=\mb{Kernel}(x)\oplus \mb{Image}(x)$ and the restriction
$x:\,\mb{Image}(x)\to \mb{Image}(x)$ is locally unipotent.\vspace*{1ex}

The {\it additive Jordan-Chevalley decomposition}: 
Let $V$ be a $\F$-linear space, and $x\in End(V)$ be locally finite. Then there exists a unique semisimple $s\in End(V)$, a unique locally 
nilpotent $n\in End(V)$ such that $x=s+n$ and $[s,n]=0$.
If $U\subseteq V$ is a $x$-invariant subspace, then $U$ is both $s$ and
$n$-invariant. Furthermore, the additive Jordan-Chevalley decomposition of the
restriction of $x$ to an endomorphism of $U$ is given by the restrictions of $s$
and $n$ to endomorphisms of $U$.\vspace*{1ex}

Let $e\in End(V)$ be an idempotent. Then
\begin{eqnarray*}
    End(V)_e:=e\, End(V)\, e=\Mklz{m\in End(V)}{me=em=m}
\end{eqnarray*}
is a monoid with unit $e$. Its unit group $(End(V)_e)^\times$ consists of the endomorphisms $m\in End(V)$ which satisfy
$\mb{Kernel}(m)=\mb{Kernel}(e)$, $\mb{Image}(m)=\mb{Image}(e)$, and for which the
restriction $ m:\mb{Image}(m)\to\mb{Image}(m)$ is invertible.
In particular $(End(V)_{id_V})^\times=Aut(V)$. 

We need the following variant of the {\it multiplicative Jordan-Chevalley decomposition}:
Let $e\in End(V)$ be idempotent, and let $m\in (End(V)_e)^\times$ be locally finite.
Then there exists a unique semisimple $s\in (End(V)_e)^\times$, a unique locally 
weak unipotent $u\in (End(V)_e)^\times$ such that $x=su =us$.\vspace*{1ex}

Let $V$, $W$ be $\F$-linear spaces. Let $e\in End(V)$, $f\in End(W)$ be idempotents, and let $x\in
(End(V)_e)^\times$, $y\in (End(W)_f)^\times$. Then $x\otimes y\in
(End(V\otimes W)_{e\otimes f})^\times$. 
If the multiplicative Jordan-Chevalley decomposition of $x$ is given by
$(s_x,\,u_x)$, and the decomposition of $y$ is given by $(s_y,\,u_y)$, then the
multiplicative Jordan-Chevalley decomposition of $x\otimes y$ is given by
$(s_x\otimes s_y,\,u_x\otimes u_y)$.
\subsection*{Some notation}
(1) If $M$ is a monoid we denote by $M^{op}$ the opposite monoid. Similarly, if $\g$ is a Lie algebra we denote by $\g^{op}$ the opposite Lie algebra. We
denote by $U(\g)^{op}$ the opposite universal enveloping algebra of $\g$. Note that $U(\g^{op})=U(\g)^{op}$.\vspace*{1ex}
 
(2) In this paper we often find the situation, where $M$ is a monoid with
coordinate ring $\FK{M}$ and associated Lie algebra $Lie(M)$, and the monoid $M^{op}\times M$ as well as the algebra $U(Lie(M))^{op}\otimes U(Lie(M))$ act on
$\FK{M}$. We denote both actions by $\pi$. If we speak of the $M$-module $\FK{M}$ we equip $\FK{M}$ with the action
\begin{eqnarray*}
  x_\ro f:=\pi(1,x)f \quad \mb{where }\quad x\in M,\;f\in\FK{M}.
\end{eqnarray*}    
If we speak of the $M^{op}$-module $\FK{M}$ we equip $\FK{M}$ with the action
\begin{eqnarray*}
  x_\lo f:=\pi(x,1)f \quad \mb{where }\quad x\in M^{op},\;f\in\FK{M}.
\end{eqnarray*}   
 
Similarly, if we speak of the $U(Lie(M))$-module $\FK{M}$, resp. of the $Lie(M)$-module $\FK{M}$, we equip $\FK{M}$ with the action
\begin{eqnarray*}
  x_\ro f:=\pi(1\otimes x)f \quad \mb{where }\quad x\in U(Lie(M)), \;\mb{resp. } x\in Lie(M),\;\mb{ and }f\in\FK{M}.
\end{eqnarray*}    
If we speak of the $U(Lie(M))^{op}$-module $\FK{M}$, resp. of the $Lie(M)^{op}$-module $\FK{M}$, we equip $\FK{M}$ with the action
\begin{eqnarray*}
  x_\lo f:=\pi(x\otimes 1)f \quad \mb{where }\quad x\in U(Lie(M)^{op}),\;\mb{resp. } x\in Lie(M)^{op},\;\mb{ and }f\in\FK{M}.
\end{eqnarray*}    
Because $U(Lie(M))^{op}\otimes U(Lie(M))$ identifies with $U(Lie(M)^{op}\times Lie(M))$, giving an action of the algebra $U(Lie(M))^{op}\otimes U(Lie(M))$ is equivalent to 
giving an action of the Lie algebra $Lie(M)^{op}\times Lie(M)$. We denote this action by $\ti{\pi}$, i.e.,
\begin{eqnarray*}
   \ti{\pi}(x,y)=\pi(x\otimes1)+\pi(1\otimes y) \quad \mb{ where }\quad x\in Lie(M)^{op},\;y\in Lie(M).
\end{eqnarray*}

We apply the same principles of notation also to other situations. For example
to the natural action of $U(\g)^{op}\otimes U(\g)$ on $U(\g)^*$, where $\g$ is
a Lie algebra.\vspace*{1ex}

(3) If $\g$ is a Lie algebra and $\mathcal C$ is a category of $\g$-modules, resp. a set of $\g$-modules 
we say $\g$ acts faithfully on the $\g$-modules contained in 
$\mathcal C$ if the zero of $\g$ is the only element which acts as zero on every $\g$-module contained in $\mathcal C$.

Similarly, if $M$ is a monoid and $\mathcal C$ is a category of $M$-modules,
resp. a set of $M$-modules we say $M$ acts faithfully on the $M$-modules contained in 
$\mathcal C$ if for all $m,\ti{m}\in M$ from
\begin{eqnarray*}
  m v=\ti{m} v \quad\mb{ for all } M\mb{-modules }V \mb{ contained in }{\mathcal C},\;v\in V
\end{eqnarray*}
follows $m=\ti{m}$.

%
\section{\label{TK} The Tannaka reconstruction}
%
%
%
In this section we describe and investigate the integration of Lie algebras by the Tannaka reconstruction. We only give the theory. 
We do not treat non-classical examples here, because this would make the paper too long. It would
also require parts of the theory of this section, and parts of the algebraic geometry of the
next section before. For some examples please look at \cite{M1} and \cite{M4}.

\subsection{The Tannaka monoid $M$, its coordinate ring $\FK{M}$, and its Lie algebra $Lie(M)$}
Let $\g$ be a Lie algebra over a field $\F$ of characteristic 0. 
Let $\mathcal C$ be a full subcategory of the category of of $\g$-modules with the following properties:
\begin{itemize}
\item[(1)] If $V$ is a $\g$-module contained in $\mathcal C$, then also every $\g$-module isomorphic to $V$ is contained in $\mathcal C$.
\item[(2)] If $V$, $W$ are $\g$-modules contained in $\mathcal C$, then $\mathcal C$ contains a sum, and a tensor product of $V$ and $W$. Furthermore, $\mathcal C$ 
contains a one-dimensional trivial $\g$-module. 
\item[(3)] If $V$ is a $\g$-module contained in $\mathcal C$, then also every $\g$-submodule of $V$ is contained in $\mathcal C$.
\item[(4)] The Lie algebra $\g$ acts faithfully on the $\g$-modules contained in $\mathcal C$, i.e., the zero of $\g$ is the only element which acts as zero on every 
$\g$-module contained in $\mathcal C$.
\end{itemize}

It is convenient to use the following notation: Let $V$ be a $\g$-module contained in $\mathcal C$. It extends to a module of the universal enveloping
algebra $U(\g)$ of $\g$. For $x\in\g$ resp. $x\in U(\g)$ we denote by $x_V$ the endomorphism of $V$ corresponding to $x$.\vspace*{1ex}
 
To associate to the category $\mathcal C$ a monoid with coordinate ring  we have to specify a category of duals in addition.
The prototypical example is the category of full duals ${\mathcal C}^{full}$ obtained as follows:
Let $V$ be a $\g$-module contained in $\mathcal C$. Its $\F$-linear dual $V^*$, which we call the {\it full dual} of $V$, 
gets the structure of a $\g^{op}$-module by the dual maps $x_V^*:V^*\to V^*$
of the maps $x_V:V\to V$, $x\in\g$.
The objects of the {\it category of full duals} ${\mathcal C}^{full}$ are the $\g^{op}$-modules $V^*$,
where $V$ is a $\g$-module contained in $\mathcal C$.
The morphisms from the object $W^*$ to the object $V^*$ are the $\g^{op}$-equivariant
linear maps given by the dual maps $\al^*:W^*\to V^*$  of the morphisms $\al:V\to W$ of $\mathcal C$.

Now we abstract the main properties of this example. A {\it category of duals} ${\mathcal C}^{du}$ is a category of $\g^{op}$-modules and
morphisms of $\g^{op}$-modules obtained in the following way:
For every object $V$ of $\mathcal C$ choose a subspace
\begin{eqnarray*}
   V^{du} \subseteq  V^*, 
\end{eqnarray*}
which separates the points of $V$, (i.e., the pairing $V^{du}\times V\to \F$ given by evaluation is non-degenerate), such that the following conditions are satisfied:
\begin{itemize}
\item[(1)] For all $x\in \g$ and all objects $V$ of $\mathcal C$ the dual map $(x_V)^{du}:V^{du}\to V^{du}$ of $x_V :V\to V$ exists. 
For all morphisms $\al:V\to W$ of $\mathcal C$ the dual map $\al^{du}:W^{du}\to V^{du}$ exists.
\item[(2)] For all $\g$-modules $V$, $W$ contained in $\mathcal C$, for all sums $V\oplus W$ and tensor products $V\otimes W$ of $V$, $W$, we have
\begin{eqnarray*}
   V^{du}\oplus W^{du} = (V\oplus W)^{du} \quad \mb{ and }\quad V^{du}\otimes W^{du} \subseteq (V\otimes W)^{du},
\end{eqnarray*}
where $V^{du}\oplus W^{du}$ resp. $V^{du}\otimes W^{du}$ are realized in the obvious way as spaces of linear functions on $V\oplus W$ resp. $V\otimes W$.
\end{itemize}
The objects of ${\mathcal C}^{du}$ are the $\g^{op}$-modules $V^{du}$, $V$ an object of $\mathcal C$. The morphisms from the object
$W^{du}$ to the object $V^{du}$ are the $\g^{op}$-equivariant linear maps given by the dual maps $\al^{du}:W^{du}\to V^{du}$  of the
morphisms $\al:V\to W$ of $\mathcal C$.\vspace*{1ex}

If $V$, $W$ are $\g$-modules contained in $\mathcal C$, then the sum of $V$, $W$ and also the tensor product of $V$, $W$ is only determined 
up to isomorphy of $\g$-modules. In this paper we only encounter situations, where a property, which holds for one sum resp. tensor product of $V$ and $W$, 
also holds for all sums resp. tensor products of $V$ and $W$. To relax our notation we do not write this down in the following definitions, theorems, and proofs. 
\vspace*{1ex}

Now fix a category of duals ${\mathcal C}^{du}$. To cut short our notation we define
\begin{eqnarray*}
      End_{V^{du}}(V):=\Mklz{\phi\in End(V)}{\mb{the dual map }\;\phi^{du}:V^{du}\to V^{du} \;\mb{ exists }     }.   
\end{eqnarray*}
Let ${\mathcal Vect}$ be the category of $\F$-linear spaces. Let $fg:{\mathcal
  C}\to {\mathcal Vect}$ be the functor forgetting the $\g$-module
structure of the objects of $\mathcal C$. Similarly, let $fg^{du}: {\mathcal C}^{du}\to {\mathcal Vect}$ be the functor forgetting the
$\g^{op}$-module structure of the objects of ${\mathcal C}^{du}$. Denote by $Nat$ the natural transformations of $fg$, which 
induce natural transformations of $fg^{du}$. Explicitely, $Nat$ consists of the families of linear maps
\begin{eqnarray*}
   m &=& \left(\,m_V\in End_{V^{du}}(V)\,\right)_{\,V \;obj.\;of\; {\mathcal C}}\;,
\end{eqnarray*} 
such that the diagram 
\begin{eqnarray*}
   V &\stackrel{m_V}{\to} & V \\
  \phi \downarrow \, &                    &  \,\downarrow \phi\\
   W &\stackrel{m_W}{\to} & W  
\end{eqnarray*}
commutes for all objects $V$, $W$, and all morphisms $\phi:V\to W$ of $\mathcal C$. 

Every $\g$-module $V$ contained in $\mathcal C$ is the union of the cyclic
submodules $U(\g)v$, $v\in V$, which are also contained in $\mathcal C$.  It is not
difficult to check that every natural transformation $m$ in $Nat$ is already
uniquely determined by its endomorphisms $m_V$ on the cyclic $\g$-modules $V$ contained in $\mathcal C$.
It follows that $Nat$ is a set. Induced by the algebras of endomorphisms  $End_{V^{du}}(V)$, 
$V$ an object of ${\mathcal C}$, the set $Nat$ gets the structure of an associative $\F$-algebra with unit, therefore also the structure of a Lie 
algebra.

Note also the following property of a natural transformation $m$ in $Nat$: Let
$V_j$, $j\in J$, be $\g$-modules contained in $\mathcal C$. If the sum
$\bigoplus_{j\in J} V_j$ is also a $\g$-module contained in $\mathcal C$ then  
\begin{eqnarray}\label{msum}
   m_{\,\bigoplus_{j\in J} V_j} = \bigoplus_{j\in J}\; m_{V_j} .
\end{eqnarray}

The {\it Tannaka monoid} $M$ introduced in the following proposition is the biggest monoid acting reasonably on the modules contained in $\mathcal C$, 
compatible with the duals in $\mathcal C^{du}$. Its coordinate ring is the {\it coordinate ring of matrix coefficients}.
\begin{Proposition}\label{TK1} Let $M$ be the set of elements $m\in Nat$ which satisfy 
\begin{itemize}
\item[(1)] $\quad m_{V\otimes W} = m_V\otimes m_W$ for all objects $V$, $W$ of $\mathcal C$,
\item[(2)] $\quad m_{V_0} = id_{V_0}$ for every trivial one-dimensional $\g$-module $V_0$.
\end{itemize}
For every object $V$ of $\mathcal C$, $v\in V$, and $\phi\in V^{du}$ define a function $f_{\phi v}: M\to\F$, the matrix coefficient of $\phi$ and $v$, by
\begin{eqnarray*}
   f_{\phi v}(m):= \phi(m_V v)\quad\mb{ where }\quad m\in M.
\end{eqnarray*}
Then $M$ is a submonoid of $Nat$ and 
\begin{eqnarray*}
  \FK{M} := \Mklz{f_{\phi v}\,}{\,\phi\in V^{du}\,,\, v\in V\,,\,\,V \;an\; object\;of\; {\mathcal C} }.
\end{eqnarray*}
is a coordinate ring on $M$. Left and right multiplications with elements of $M$ induce comorphisms of $\FK{M}$. As a consequence the monoid $M^{op}\times M$ acts on 
the algebra $\FK{M}$ by 
\begin{eqnarray*}
   \left(\pi(m_1,m_2)f\right)(m):=f(m_1mm_2)\quad\mb{ where } \quad f\in\FK{M},\;\mb { and }\; m_1,\,m_2,\,m\in M. 
\end{eqnarray*}
\end{Proposition}
\begin{Remark} Assigning to an element of $M$ its evaluation homomorphism, $M$ embeds as a set into the $\F$-valued points of $\FK{M}$. In general, this 
embedding is not surjective. In general, the multiplication map of $M$ does not induce a comultiplication of the coordinate ring $\FK{M}$.
\end{Remark}

\Proof It is trivial to check that $M$ is a submonoid of $Nat$. 
To show that hat $\FK{M}$ is an algebra of functions on $M$ let $V$, $W$ be objects of $\mathcal C$, and $v\in V$, $w\in W$, $\phi\in
V^{du}$, and $\psi\in W^{du}$. Because of property (1) of the definition of $M$ we have $f_{\phi v}f_{\psi w} = f_{\phi\otimes\psi\, v\otimes w} $.
Because of property (2) we have $f_{\phi v} = \phi(v)\, 1$ for all $\phi\in V_0^{du} =V_0^*$, $v\in V_0$.

Let $m,\,\ti{m}\in M$, and suppose $f(m)=f(\ti{m})$ for all $f\in\FK{M}$. Then for every object $V$ of $\mathcal C$ we have 
$\phi(m_V v)=\phi(\ti{m}_V v)$ for all $\phi\in V^{du}$, $v\in V$. Because $V^{du}$ is point separating on V, we find $m_V=\ti{m}_V$.

Let $m\in M$. The left multiplication $l_m$ and the right multiplication $r_m$ induce comorphisms, because for $v\in V$, $\phi\in V^{du}$, 
and $V$ an object of $\mathcal C$ we have $l_m^* \,f_{\phi v} = f_{m_V^{du}(\phi)\, v}$ and  $r_m^* \,f_{\phi v}=f_{\phi\, m_V(v)}$.
\qed

Let $V$ be an object of $\mathcal C$. Equip $V$ with the coordinate ring $\FK{V}$, which is generated by the functions of $V^{du}$. In the obvious way
$\FK{V}$ is isomorphic to the symmetric algebra in $V^{du}$. The following proposition is trivial to show:

\begin{Proposition}\label{TK2} The monoid $M$ acts on $V$, such that the for all $m\in M$ the left applications $l_m:V\to V$ defined by 
$l_m(v):=m_V v$, $v\in V$, and for all $v\in V$ the evaluations $r_v:M\to V$ defined by $r_v(m):=m_V v$,  $m\in M$, 
are morphisms of sets with coordinate rings.
\end{Proposition}

The {\it Lie algebra} $Lie(M)$ of ($M$, $\FK{M}$) is the biggest Lie algebra acting reasonably on the 
modules of $\mathcal C$, compatible with the duals in ${\mathcal C}^{du}$, and with the coordinate ring $\FK{M}$.
It is introduced in the following proposition together with the {\it adjoint action} of the unit group $M^\times$. 

\begin{Proposition}\label{TK3} Let $Lie(M)$ be the set of natural transformations $x\in Nat$, which satisfy the following properties:
\begin{itemize}
\item[(1)] $\quad x_{V\otimes W}= x_V\otimes id_{W}+ id_V \otimes x_W$ for all objects $V$, $W$.
\item[(2)] $\quad x_{V_0} = 0_{V_0}$ for every trivial one-dimensional $\g$-module $V_0$.
\item[(3)] There exists a map $\delta_x:\FK{M}\to \F$ such that 
\begin{eqnarray*}
 \qquad \qquad\delta_x(f_{\phi v})=\phi(x_V v)\quad \mb{ for all } \quad \phi\in V^{du},\;v\in V, \;V \mb{ an object of } {\mathcal C}.
\end{eqnarray*}
\end{itemize}
Then for every $x\in Lie(M)$ the map $\delta_x$ is a derivation of $\FK{M}$ in
$1\in M$. Furthermore, $Lie(M)$ is a Lie subalgebra of the Lie algebra $Nat$. The unit group $M^\times$ of $M$ acts on $Lie(M)$ by conjugation.
\end{Proposition}
\begin{Remark} Assigning to an element $x\in Lie(M)$ the derivation $\delta_x\in Der_1(\FK{M})$, the linear space $Lie(M)$ embeds into the linear space 
$Der_1(\FK{M})$. In general, this embedding is not surjective.
\end{Remark}

\Proof It is not difficult to see that for an element $x\in Nat$ the properties (1) and (2) and (3) are equivalent to any of the following properties:
\begin{itemize}
\item[(4)] There exists a derivation $\delta_x:\FK{M}\to \F$ in $1\in M$ such that 
\begin{eqnarray*}
 \qquad \qquad\delta_x(f_{\phi v})=\phi(x_V v)\quad\mb{ for all } \quad \phi\in V^{du},\;v\in V, \;V \mb{ an object of } {\mathcal C}.
\end{eqnarray*}
\item[(4')] There exists a derivation $\pi(1\otimes x)\in Der(\FK{M})$ such that
\begin{eqnarray*} 
   \qquad\qquad \pi(1\otimes x) (f_{\phi v}) = f_{\phi\, x_V v} \quad \mb{ for all } \quad\phi\in V^{du},\; v\in V, \;V \mb{ an object of } {\mathcal C}\;.
\end{eqnarray*}
\end{itemize}
\begin{itemize}
\item[(4'')] There exists a derivation $\pi(x\otimes 1)\in Der(\FK{M})$ such that
\begin{eqnarray*} 
   \qquad\qquad \pi(x\otimes 1)(f_{\phi v}) = f_{x_V^{du}\phi\,v}\quad \mb{ for all } \quad\phi\in V^{du},\; v\in V, \;V \mb{ an object of } {\mathcal C}\;.
\end{eqnarray*}
\end{itemize}

If $x,\,y\in Lie(M)$ then also $[x,y]\in Lie(M)$, since the derivation $[\pi(1\otimes x),\pi(1\otimes y)]$ satisfies property (4') for $[x,y]$.
Let $x\in Lie(M)$ and $m\in M^\times$. The comorphisms $l_m^*$, $r_{m^{-1}}^*$ of the left and right multiplications $l_m$, $r_{m^{-1}}$ by $m$, 
$m^{-1}$ are isomorphism of algebras. Therefore $\delta_x\circ l_m^*\circ r_{m^{-1}}^* $ satisfies property of (4) for $mxm^{-1}$.
\qed

The Lie algebra $Lie(M)$ acts on every object $V$ of $\mathcal C$. Later we will interprete this action as the differentiated action of the monoid
$M$ on $V$. 
With $Lie(M)$ also the universal enveloping algebra $U(Lie(M))$ acts on the objects of $\mathcal C$. 
For $x\in U(Lie(M))$ we denote by $x_V$ the corresponding endomorphism in
$End_{V^{du}}(V)$, $V$ an object of $\mathcal C$. From the proof of the last proposition follows immediately:

\begin{Corollary} We get an action $\pi$ of $U(Lie(M))^{op}\otimes U(Lie(M))$ on $\FK{M}$ by
\begin{eqnarray*}
         \pi(x\otimes y) f_{\phi v} = f_{x_V^{du}\phi\, y_V v} 
\end{eqnarray*}
where  $x\in U(Lie(M))^{op}$, $y\in U(Lie(M))$, and $\phi\in V^{du}$, $v\in V$, $V$ an object of $\mathcal C$.

Recall our notation for such an action from Section \ref{Pre}. The Lie algebra $Lie(M)$ acts on $\FK{M}$ by left invariant derivations, 
i.e., the derivations commute with the homomorphisms given by the action of $M^{op}$.  
The Lie algebra $Lie(M)^{op}$ acts on $\FK{M}$ by right invariant derivations, 
i.e., the derivations commute with the homomorphisms given by the action of $M$.
\end{Corollary}

Let $N$ be a submonoid of $M$. Denote by $I(N)$ the vanishing ideal of $N$ in
$\FK{M}$. It is easy to check that
\begin{eqnarray*}
   Lie(N):=\Mklz{x\in Lie(M)}{ \delta_x(I(N))=0}
\end{eqnarray*}
is a Lie subalgebra of $Lie(M)$, which we call the {\it Lie algebra} of $N$.
Note that $Lie(N)$ consists of the elements of $x\in Lie(M)$, for which the corresponding
derivation $\delta_x$ of $\FK{M}$ in $1\in M$ factors to a derivation of $\FK{N}$ in
$1\in N$. In this way, we get an embedding of the linear space $Lie(M)$ into $Der_1(\FK{N})$.\vspace*{1ex}

The monoid $M$ and its coordinate ring of matrix coefficients $\FK{M}$ may be quite small. The next definition singles out when we can look at 
$(M,\,\FK{M})$ as an analogue of a linear algebraic monoid integrating $\g$.
By our assumptions on the category $\mathcal C$ the map
\begin{eqnarray*}
   \g &\to& Nat\\
    x &\mapsto & (x_V)_{\,V \;obj.\;of\; {\mathcal C}}
\end{eqnarray*}
is an embedding of Lie algebras. {\bf Please note:} In the following we identify $\g$ with its
image in $Nat$ without further mentioning. 
The elements $\g\subseteq Nat $ satisfy properties (1) and (2) of Proposition \ref{TK3}, but in general not property (3).
\begin{Definition} We call the pair $\mathcal C$, ${\mathcal C}^{du}$ good for integrating $\g$, if $\g\subseteq Lie(M)$.
We call the pair $\mathcal C$, ${\mathcal C}^{du}$ very good for integrating $\g$, if $\g\subseteq Lie(M^\times)$.
\end{Definition}

To realize $\g$ as a subalgebra of the Lie algebra of the unit group is well known from linear algebraic monoids. 
In all examples which we treat in \cite{M1} and \cite{M4} the pair $\mathcal C$, ${\mathcal C}^{du}$ is very good for integrating $\g$. 
Nevertheless there may exist other examples where $\mathcal C$, ${\mathcal C}^{du}$ is only good for integrating $\g$.
Nearly all theorems which we treat in this paper only assume that the pair $\mathcal C$, ${\mathcal C}^{du}$ is good for integrating $\g$.
The relation between the two definitions is described by Proposition
\ref{Defrel}. For its proof we first note the following easy consideration
which will be useful several times in this paper. 
\begin{Proposition} \label{sub} Let $V$ be a $\g$-module contained in
  $\mathcal C$. Every $\g$-invariant subspace $U$ of $V$ is also $M$-invariant. 
\end{Proposition}

\Proof Denote by $j:\,U\to V$ be the inclusion map. By our assumptions on the category $\mathcal C$ the $\g$-module $U$ is contained in $\mathcal C$. 
Therefore, it gets the structure of a $M$-module. The $\g$-equivariant morphism $j:\,U\to V$ is also a $M$-equivariant morphism. 
Therefore, $U=j(U)$ is a $M$-invariant submodule of $V$.\qed

\begin{Proposition}\label{Defrel} It is equivalent:
\begin{itemize}
\item[(i)] The pair $\mathcal C$, ${\mathcal C}^{du}$ is very good for integrating $\g$.
\item[(ii)] The pair $\mathcal C$, ${\mathcal C}^{du}$ is good for integrating $\g$ and $Lie(M^\times)=Lie(M)$.
\item[(iii)] The pair $\mathcal C$, ${\mathcal C}^{du}$ is good for integrating $\g$ and $M^\times$ is dense in $M$.
\end{itemize}
\end{Proposition}

\Proof Obviously from (iii) follows (ii). From (ii) follows (i) by
definition. Now suppose that (i) holds. To show (iii) we have to show that $M^\times$ is
dense in $M$. 
Let $V$ be a $\g$-module contained in $\mathcal C$, let $v\in V$ and $\phi\in V^{du}$ such that $f_{\phi  v}(M^\times)=\{0\}$.
Because of $\g\subseteq Lie(M^\times)$, the vanishing ideal $I(M^\times)$ is stable under the derivations $\pi(1\otimes x)$, $x\in\g$. Therefore
$I(M^\times)$ is also stable under $\pi(1\otimes x)$, $x\in U(\g)$. It follows 
\begin{eqnarray*}
    f_{\phi \,x_V v}=\pi(1\otimes x)f_{\phi v}\in I(M^\times) \quad \mb{ for all } \quad x\in U(\g).
\end{eqnarray*}
By evaluating in $1\in M^\times$ we find $\phi(U(\g)v)=\{0\}$. By Proposition \ref{sub} the
$\g$-invariant subspace $U(\g)v$ of $V$ is also $M$-invariant. It follows $f_{\phi v}(M)=\phi(Mv)=\{0\}$.
\qed

Note also the following saturation property, which can be used sometimes to
replace in a statement or conclusion the Lie algebra $\g$ by $Lie(M)$: Let the pair $\mathcal C$, ${\mathcal C}^{du}$ be good for integrating $\g$. Then every $\g$-module 
contained in $\mathcal C$ extends to a $Lie(M)$-module. The $\g$-equivariant maps between $\g$-modules contained in $\mathcal C$
coincide with the $Lie(M)$-equivariant maps. In this way the category $\mathcal C$ gives a category of $Lie(M)$-modules. Similarly, the category 
${\mathcal C}^{du}$ gives a category of $Lie(M)^{op}$-modules. We can now use these new categories for the Tannaka
reconstruction. We get back the monoid $M$ and its coordinate ring $\FK{M}$. The associated Lie algebra is $Lie(M)$.

\subsection{The kernel of the adjoint action of $M^\times$ on $Lie(M)$}
Denote by $C(M)$ resp. $C(M^\times)$ the center of $M$ resp. $M^\times$.
\begin{Proposition} The center $C(M)$ is a closed abelian submonoid of $M$.
\end{Proposition}

\Proof Because of $C(M) =\bigcap_{n\in M}\Mklz{m\in M}{mn=nm }$
it is sufficient to show that for any $n\in M$ the set $\Mklz{m\in M}{mn=nm }$ is closed. Since the duals are point separating it 
coincides with the common zero set of the functions
\begin{eqnarray*}
   f_{\phi\,n_V v}-f_{n_V^{du}\phi\,v}   \quad\mb{ where }\quad \phi\in V^{du},\;v\in V,\;V\mb{ an object of }{\mathcal C}.
\end{eqnarray*}
\qed
\begin{Theorem} If $\mathcal C$, ${\mathcal C}^{du}$ is good for integrating $\g$ then $C(M)\cap M^\times$ is the kernel of the adjoint 
action of $M^\times$ on $Lie(M)$. 
If $\mathcal C$, ${\mathcal C}^{du}$ is very good for integrating $\g$ then this kernel coincides with $C(M^\times)$. 
\end{Theorem}

\Proof (a) We first show
\begin{eqnarray}\label{C(M)Lie}
   C(M) =\Mklz{m\in M}{mx=xm \mb{ for all } x\in Lie(M) } .
\end{eqnarray}
Let $m\in C(M)$. Then $m\ti{m}=\ti{m}m$ for all $\ti{m}\in M$, from which follows
\begin{eqnarray*}
   \phi(m\ti{m}v)=\phi(\ti{m}m v) \quad\mb{ for all } \quad \ti{m}\in M,\;\phi\in V^{du},\;v\in V,\;V\mb{ an object of }{\mathcal C}. 
\end{eqnarray*}
This is equivalent to
\begin{eqnarray*}
   f_{m_V^{du}\phi \,v }=f_{\phi\,m_V v} \quad\mb{ for all } \quad \phi\in V^{du},\;v\in V,\;V\mb{ an object of }{\mathcal C}. 
\end{eqnarray*}
Now let $x\in Lie(M)$. By applying the derivation $\delta_x$ to this equation we get
\begin{eqnarray*}
    \left(m_V^{du}\phi\right) (x_V v )= \phi(x_V m_V v) \quad\mb{ for all } \quad \phi\in V^{du},\;v\in V,\;V\mb{ an object of }{\mathcal C}. 
\end{eqnarray*}
Since the duals are point separating it follows $mx=xm$.\vspace*{1ex}

Now let $m\in M$ such that $mx=xm$ for all $x\in Lie(M)$. By assumption
$\g\subseteq Lie(M)$. Therefore, for every object $V$ contained in $\mathcal
C$, the map from $V$ to $V$, which maps the element $v\in V$ to $m_V v\in V$, 
is a morphism of $\g$-modules. Because of $M\subseteq Nat$ we get $m_V\ti{m}_V=\ti{m}_V m_V$ for all $\ti{m}\in M$. Therefore 
$m\ti{m}=\ti{m}m$ for all $\ti{m}\in M$.\vspace*{1ex}

(b) From equality (\ref{C(M)Lie}) follows that $C(M)\cap M^\times$ is the kernel of the adjoint action. Obviously 
$C(M)\cap M^\times\subseteq C(M^\times)$. Now let $\mathcal C$, ${\mathcal
  C}^{du}$ be very good for integrating $\g$. If $m\in C(M^\times)$ then 
$ m\ti{m}=\ti{m}m$ for all $\ti{m}\in M^\times$, from which follows
\begin{eqnarray*}
  f_{m_V^{du}\phi\,v}(\ti{m})=f_{\phi\,m_V v} (\ti{m})\quad 
                    \mb{ for all } \quad\ti{m}\in M^\times,\;\phi\in V^{du},\;v\in V,\;V\mb{ an object of }{\mathcal C}. 
\end{eqnarray*}
By Proposition \ref{Defrel} the unit group
$M^\times$ is dense in $M$. Therefore these equations are also valid for all $\ti{m}\in M$. Since the duals are point separating it follows  
$ m\ti{m}=\ti{m}m $ for all $\ti{m}\in M$.
\qed
\subsection{The relation between $M$ and $\g$-invariant subspaces}
For one direction recall Proposition \ref{sub}. The next theorem shows the
other direction for closed subspaces.

\begin{Theorem}\label{sub2} Let $\mathcal C$, ${\mathcal C}^{du}$ be good for integrating
  $\g$. Let $V$ be a $\g$-module contained in $\mathcal C$. Then every 
  $M$-invariant closed subspace $U$ of $V$ is also $\g$-invariant.
\end{Theorem}
\begin{Remark}\label{RemInv} If $V^{du}=V^*$, in particular if $V$ is finite dimensional, then every subspace of $V$ is closed. In this
  case the $\g$-invariant subspaces of $V$ and the $M$-invariant subspaces of $V$ coincide.
  Later we introduce the class of examples where $\g$ is generated by integrable locally finite
  elements. Also in this case the $\g$-invariant subspaces of $V$ and the $M$-invariant subspaces of
  $V$ coincide. Compare Theorem \ref{sub3}.
\end{Remark}

\Proof Let $x\in\g$. Let $h\in \FK{V}$ such that $h(u)=0$ for all $u\in U$. We
show that then also $h(xu)=0$ for all $u\in U$. Since $U$ is closed this
implies $xU\subseteq U$.

The algebra $\FK{V}$ is a symmetric algebra in $V^{du}$. Write $h$ in the form
$h=\sum_{n=0}^N h_n$ where $h_n\in\FK{V}$ is homogeneous in degree
$n$. Because of
\begin{eqnarray*}
  0=h(tu)=\sum_{n=0}^N t^n h_n(u) \quad\mb{ for all }\quad u\in U,\;t\in\F,
\end{eqnarray*}
we find $h_n(u)=0$ for all $u\in U$, $n=0,\,1,\,\ldots,\, N$. 
Therefore it is sufficient to assume that $h\in \FK{V}$ is homogeneous of
degree $n\in\N$, i.e., 
\begin{eqnarray*}
  h=\sum_{i=1}^m \phi_1^{(i)}\phi_2^{(i)}\cdots\phi_n^{(i)}\quad \mb{ with }\quad
  \phi_j^{(i)}\in V^{du},
\end{eqnarray*}
such that $h(u)=0$ for all $u\in U$. Set 
\begin{eqnarray*}
  \ti{h}(u_1,\,u_2,\,\ldots,\,u_n):=\sum_{i=1}^m
  \phi_1^{(i)}(u_1)\phi_2^{(i)}(u_2)\cdots\phi_n^{(i)}(u_n) \quad \mb{ for
  }\quad u_1,\,\ldots,\,u_n\in U.
\end{eqnarray*}
Denote by $S_n$ the permutation group in $n$ variables. We have
\begin{eqnarray*}
  0=h(t_1 u_1 + t_2 u_2+ \cdots t_n u_n) \quad\mb{ for all }\quad t_1,\,\ldots,\,t_n\in \F,\;\; u_1,\,\ldots,\,u_n\in U.
\end{eqnarray*}
The expression on the right is polynomial in $t_1,\,\ldots,\,t_n$. Comparing
the coefficients of $t_1 t_2\cdots t_n$ we get
\begin{eqnarray}\label{PerF}
  0=\sum_{\sigma\in S_n} \ti{h}(u_{\sigma 1},\,\ldots,\,u_{\sigma n} )\quad
  \mb{ for all }\quad u_1,\,\ldots,\,u_n\in U.
\end{eqnarray}
Since $U$ is
$M$-invariant we can replace in this equation $u_1$ by $m u_1$, $m\in M$. In
this way we get a function of $\FK{M}$ which is zero. Therefore also the
derivation $\delta_x$ of $\FK{M}$ in $1\in M$ applied to this function is
zero. This means that we can replace in equation (\ref{PerF}) $u_1$ by
$xu_1$. Repeating the same argument we can replace in equation (\ref{PerF}) $u_2$ by
$xu_2$, $\ldots$, $u_n$ by $xu_n$. Setting
$u_1=u_2=\ldots=u_n=:u$ we get
\begin{eqnarray*}
  0=n!\,\ti{h}(xu,\,\ldots,\, xu)=n!\,h(xu) \quad\mb{ for all }\quad u\in U.
\end{eqnarray*}
\qed
\subsection{The embedding of $\FK{M}$ into the algebra $U(\g)^*$}
Recall from Section \ref{Pre} the canonical Hopf algebra structure of the
universal enveloping algebra $U(\g)$ of $g$. We denote by $\Delta:\,U(\g)\to U(\g)\otimes U(\g)$ 
the comultiplication, by $\eps:\,U(\g)\to\F$ the counit, and by $S:\,U(\g)\to
U(\g)$ the antipode. Recall also from Section \ref{Pre} the commutative
algebra structure of the full dual $U(\g)^*$.
We get an action of $U(\g)^{op}\otimes U(\g)$ on $U(\g)^*$ by
\begin{eqnarray*}
           \left( \pi(a\otimes b) h\right) (x) := h(axb)\quad\mb{ where }\quad a,\,b\in\g,\;h\in U(\g)^*,\;x\in U(\g).
\end{eqnarray*} 
Recall our notation for this action from Section \ref{Pre}. The Lie algebra $\g$ acts on $U(\g)^*$
by derivations. Similarly, the Lie algebra $\g^{op}$ acts on $U(\g)^*$ by derivations.
\begin{Theorem}\label{UnivEin}
(a) For $V$ an object of $\mathcal C$, $v\in V$, and $\phi\in V^{du}$ define a function $g_{\phi v}: U(\g)\to\F$, the matrix coefficient of 
$\phi$ and $v$, by
\begin{eqnarray*}
   g_{\phi v}(x):= \phi(x_V v)\quad\mb{ where }\quad x\in U(\g).
\end{eqnarray*}
We get a $U(\g)^{op}\otimes U(\g)$-invariant subalgebra of $U(\g)^*$, the
algebra of matrix coefficients on $U(\g)$, by 
\begin{eqnarray*}
    \FK{U(\g)}:=\Mklz{g_{\phi v}}{\phi\in V^{du},\;v\in V,\;V\mb{ an object of }\mathcal C}.
\end{eqnarray*}

(b) Let $\mathcal C$, ${\mathcal C}^{du}$ be good for integrating $\g$. Restrict the 
$U(Lie(M))^{op}\otimes U(Lie(M))$-action on the coordinate ring $\FK{M}$ to a $U(\g)^{op}\otimes U(\g)$-action. We get a
$U(\g)^{op}\otimes U(\g)$-equivariant isomorphism of algebras
\begin{eqnarray*}
       \quad   \Psi:\,\FK{M}\to \FK{U(\g)}
\end{eqnarray*} 
by
\begin{eqnarray*}
  \Psi(f)(x):=(x_\ro f)(1)=(x_\lo f)(1)\quad \mb{ where }\quad x\in U(\g),\;f\in\FK{M}.
\end{eqnarray*}
In particular, $\Psi(f_{\phi v})=g_{\phi v}$, $\phi\in V^{du}$, $v\in V$, $V$ an object of $\mathcal C$. 
\end{Theorem}

\Proof To (a): $\FK{U(\g)}$ is a $\F$-linear subspace of $U(\g)^*$ which is $U(\g)^{op}\otimes U(\g)$-invariant because 
$\pi(a,b)g_{\phi v}=g_{a_v^{du}\phi \,b_V v}$ for all $a\in U(\g)^{op}$, $b\in U(\g)$, and $v\in V$, $\phi\in V^{du}$, $V$ an object of $\mathcal C$.
The unit $\eps$ of $U(\g)^*$ is contained in $\FK{U(\g)}$ since
\begin{eqnarray}\label{Uunit}
  g_{\phi_0 v_0}=\phi_0(v_0)\eps
\end{eqnarray}
for all $\phi_0\in V_0^{du}$, $v_0\in V_0$, $V_0$ a one-dimensional trivial
$\g$-module. To show that $\FK{U(\g)}$ is a subalgebra of $U(\g)^*$, it is 
sufficient to show
\begin{eqnarray}\label{Umult}
  g_{\phi v}\cdot g_{\psi w}=g_{\phi\otimes \psi\, v\otimes w}
\end{eqnarray}
for all $v\in V$, $w\in W$, $\phi\in V^{du}$, $\psi\in W^{du}$, and all objects $V$, $W$ of $\mathcal C$.
We have
\begin{eqnarray*}
   \left(g_{\phi v}\cdot g_{\psi w}\right)(1)=
   (g_{\phi v}\otimes g_{\psi w})(1\otimes 1)=\phi(v)\psi(w)=g_{\phi\otimes \psi\, v\otimes w}(1).
\end{eqnarray*}
Now let $a\in U(\g)$ such that 
\begin{eqnarray*}
   \left( g_{\phi v}\cdot g_{\psi w}\right)(a)=
   (g_{\phi\otimes \psi\, v\otimes w})(a)
\end{eqnarray*}
for all $v\in V$, $w\in W$, $\phi\in V^{du}$, $\psi\in W^{du}$, and all objects $V$, $W$ of $\mathcal C$.
The Lie algebra $\g$ acts by derivations on $U(\g)^*$. Therefore, for all
$x\in\g$, for all $v\in V$, $w\in W$, $\phi\in V^{du}$, $\psi\in W^{du}$, and
all objects $V$, $W$ of $\mathcal C$ we find
\begin{eqnarray*}
  && \left(g_{\phi v}\cdot g_{\psi w}\right)(ax)
   =\left(x_\ro(g_{\phi v}\cdot g_{\psi w}) \right)(a)
   =\left((x_\ro g_{\phi v})\cdot g_{\psi w}  +g_{\phi v}\cdot (x_\ro g_{\psi w})\right)(a)\\
  &&=\left(g_{\phi\, x v}\cdot g_{\psi w}  +g_{\phi v}\cdot g_{\psi\, xw}\right)(a)
   = g_{\phi \otimes \psi \,(x v)\otimes w}(a) + g_{\phi\otimes\psi\, v\otimes (xw)}(a)
   \\&&= g_{\phi \otimes \psi \,x(v\otimes w)}(a)=  g_{\phi \otimes \psi \,v\otimes w}(ax).
\end{eqnarray*}

To (b): Obviously, $\Psi(f_{\phi v})=g_{\phi v}$ for all $\phi\in V^{du}$,
$v\in V$, $V$ an object of $\mathcal C$. It follows that $\Psi$ is a
surjective equivariant linear map. The equations (\ref{Uunit}) and (\ref{Umult}) show that it is also a morphism of algebras.
To show that $\Psi$ is injective let $V$ be an object of ${\mathcal C}$, $v\in V$, and $\phi\in V^{du}$ such that 
$\Psi(f_{\phi v})(x)=\phi(xv)=0$ for all $x\in U(\g)$. By 
Proposition \ref{sub} the $\g$-invariant subspace $U(\g)v$ is also $M$-invariant. Therefore, $Mv\subseteq U(\g)v$, from which follows 
$f_{\phi v}(m)=\phi(mv)=0$ for all $m\in M$.
\qed 

By Theorem \ref{zd} the algebra $U(\g)^*$ has no zero divisors. By the last
theorem we obtain:

\begin{Corollary} Let $\mathcal C$, ${\mathcal C}^{du}$ be good for integrating $\g$. Then $M$ is irreducible.
\end{Corollary}
\subsection{A Peter-Weyl type theorem for $\FK{M}$ if the category $\mathcal C$ is semisimple}
We call the category $\mathcal C$ {\it semisimple} if every $\g$-module belonging to $\mathcal C$ is a (possibly infinite) sum of 
irreducible $\g$-modules, which then also belong to $\mathcal C$.

Let $V$ be an irreducible $\g$-module contained in $\mathcal C$. Then $End_{\bf g}(V)$ is a skew-field, extending the field $\F$, which is identified 
with $\F id_V$. In the obvious way $V$ is a left $End_{\bf g}(V)$-module. By our assumptions on the category of duals,
$V^{du}$ gets the structure of a right $End_{\bf g}(V)$-module by 
\begin{eqnarray*}
    \phi\al:=\al^{du}\phi\quad \mb{ where }\quad\phi\in V^{du},\;\al\in End_{\bf g}(V). 
\end{eqnarray*}
To cut short our notation set $S(V):=End_{\bf g}(V)$. It is easy to check that the $\F$-linear space $V^{du}\otimes_{S(V)} V$ is in the obvious 
way a $M^{op}\times M$-module as well as a $U(Lie(M))^{op}\otimes U(Lie(M))$-module.

If the category $\mathcal C$ is semisimple this does not imply that $\mathcal
C^{du}$ is semisimple. Also, if $V$ is an irreducible $\g$-module contained in
$\mathcal C$ then this does not imply that $V^{du}$ is an irreducible
$\g^{op}$-module. Nevertheless, the following Peter-and-Weyl-type theorem holds.
\begin{Theorem}\label{PW} Let $\mathcal C$, ${\mathcal C}^{du}$ be good for integrating $\g$. 
Let $\mathcal C$ be semisimple. Let $Irr$ be a complete set of pairwise non-isomorphic irreducible modules, i.e., every irreducible module 
contained in $\mathcal C$ is isomorphic to exactly one module contained in $Irr$. Then the map 
\begin{eqnarray*}
   \Phi:\;\bigoplus_{V\in Irr} V^{du}\otimes_{S(V)} V \to \FK{M}
\end{eqnarray*}
induced by 
\begin{eqnarray*}
    \Phi(\phi\otimes v):= f_{\phi v}\quad \mb{ where } \quad \phi\in V^{du},\;v\in V,\;V\mb{ an object of }{\mathcal C},
\end{eqnarray*}
is a $M^{op}\times M$-equivariant and $U(Lie(M))^{op}\otimes U(Lie(M))$-equivariant $\F$-linear bijective map.
\end{Theorem}

\Proof It is easy to check that $\Phi$ is a well defined $M^{op}\times
M$-equivariant and $U(Lie(M))^{op}\otimes U(Lie(M))$-equivariant $\F$-linear map. 

To show that $\Phi$ is surjective consider a matrix coefficient $f_{\phi w}$, where $\phi\in W^{du}$,  $w\in W$, and $W$ an object of $\mathcal C$. Since
$\mathcal C$ is semisimple we can find a decomposition $W=\bigoplus_{i\in I}V_i$ where $V_i$ is a irreducible $\g$-module contained in $\mathcal
C$, $i\in I$. Write $w$ in the form $w=v_{i_1}+\cdots+ v_{i_n}$ with $v_{i_1}\in V_{i_1}$, $\ldots$, $v_{i_n}\in V_{i_n}$. Then
\begin{eqnarray*}
 f_{\phi w}= f_{\phi v_{i_1}}+\cdots+f_{\phi v_{i_n}} .
\end{eqnarray*}
By our assumptions on the categories $\mathcal C$ and $\mathcal
C^{du}$ the restriction of $W^{du}$ onto $V_i$ is contained in
$(V_i)^{du}$. Therefore the matrix coefficients $f_{\phi v_{i_1}}$, $\ldots$,
$f_{\phi v_{i_n}}$, and also their sum $f_{\phi w}$ are contained in the image of $\Phi$. 

Next we show $\FK{M}=\bigoplus_{V\in Irr}\Phi(V^{du}\otimes_{S(V)} V)$. For
this note first that every function of $\Phi(V^{du}\otimes_{S(V)} V)$ is of the form
\begin{eqnarray*}
  \quad \Phi(\sum_{i=1}^n\phi_i\otimes v_i) = f_{\phi_1\oplus \cdots\oplus \phi_n\,v_1\oplus\cdots \oplus v_n} \quad\mb{ with } \quad
  \phi_1,\ldots,\,\phi_n\in V^{du},  \;v_1,\,\ldots,\,v_n\in V,  \;n\in \N.
\end{eqnarray*} 

Let $V_1,\,\ldots,\, V_n\in Irr$ be pairwise non-isomorphic. Let $f_i\in \Phi((V_i)^{du}\otimes_{S(V_i)} V_i)$, $i=1,\,\ldots,\, n$, such that
$\sum_{i=1}^n f_i=0$. Then $f_i=f_{\psi_i w_i}$ for some $\psi_i\in (V_i^{\oplus m_i})^{du}$, $w_i\in (V_i)^{\oplus m_i}$, $i=1,\,\ldots,\,n$. 
The map
\begin{eqnarray*}
\Gamma:\,(V_1)^{\oplus m_1}\oplus\cdots\oplus (V_n)^{\oplus m_n} &\to & \;\;\qquad \FK{M}\\
w_1\oplus \cdots \oplus w_n \qquad\quad &\mapsto & f_{\psi_1 w_1}+\cdots+f_{\psi_m w_m}
\end{eqnarray*}
is a $\g$-equivariant $\F$-linear map. The image of $\Gamma$ is again completely reducible. The function $f_i$ is contained in the $V_i$-isotypical
component of the image of $\Gamma$. Therefore, from $\sum_{i=1}^n f_i=0$ follows $f_i=0$, $i=1,\,\ldots,\,n$.

To prove the injectivity of $\Phi$ it is now sufficient to prove the injectivity of
$\Phi$ restricted to $V^{du}\otimes_{S(V)} V$ for every $V\in Irr$. 
Fix $V\in Irr$. Let $\Psi:\FK{M}\to\FK{U(\g)}$ be the isomorphism described in Theorem \ref{UnivEin} (b). It is sufficient to show: Let $n\in \N$. 
Let $v_1,\,\ldots,\,v_n\in V$ be $S(V)$-linearly independent, let $\phi_1,\,\ldots,\,\phi_n\in V^{du}$ such that
\begin{eqnarray}\label{fmu}
    \Psi\left( f_{\phi_1 v_1} + \cdots + f_{\phi_n v_n}\right) (U(\g))= \left( g_{\phi_1 v_1} + \cdots + g_{\phi_n v_n}\right) (U(\g))=\{0\}.
\end{eqnarray}
Then $\phi_1=\cdots=\phi_n=0$.

We show this statement by induction over $n\in \N$. For $n=1$ the subspace $U(\g)v_1$ of $V$ is $\g$-invariant and non-zero. Since
$V$ is irreducible $U(\g)v_1=V$. From (\ref{fmu}) follows $\phi_1=0$.

Now let $n\in\N\setminus\{1\}$. First assume that there exists an $x\in U(\g)$ such that $x v_1,\,\ldots,\, xv_n$ are $S(V)$-linearly dependent, and at least one of 
these elements is non-zero. 

We may assume that $x v_1,\,\ldots,\, xv_k$ is a $S(V)$-linearly independent maximal subsystem of $x
 v_1,\,\ldots,\,xv_n$, $k\in\{1,\,\ldots,\, n-1\}$. Then
\begin{eqnarray*}
 xv_{j}= c_{j1}\,xv_1+\cdots + c_{jk}\,xv_k \quad \mb{ for some }\quad c_{j1},\,\ldots,\,c_{jk}\in S(V), \quad j=k+1,\,\ldots,\,n.
\end{eqnarray*} 
If we define
\begin{eqnarray*}
   \ti{v}_j:=v_j- c_{j1}\,v_1- \cdots - c_{jk} \,v_k \quad \mb{ for } \quad j=k+1,\,\ldots,\,n,\\
   \ti{\phi}_l:= \phi_l+ c_{k+1\,l}\,\phi_{k+1}+\cdots+ c_{n l}\,\phi_n\quad \mb{ for } \quad l=1,\,\ldots,\,k,\\
\end{eqnarray*} 
then also $\ti{v}_{k+1},\,\ldots,\,\ti{v}_n$ are
$S(V)$-linearly independent, $x\ti{v}_{k+1}=\cdots =x\ti{v}_n=0$, and
\begin{eqnarray}\label{mbeq}
  g_{\phi_1 v_1} + \cdots + g_{\phi_n v_n}=g_{\ti{\phi}_1 v_1} + \cdots
  g_{\ti{\phi}_k v_k}+ g_{\phi_{k+1} \ti{v}_{k+1}}+ \cdots + g_{\phi_n \ti{v}_n}.
\end{eqnarray}
With (\ref{fmu}) follows
\begin{eqnarray*}
   \{0\}= \left(g_{\phi_1 v_1} + \cdots + g_{\phi_n v_n}\right)(xU(\g)) = \left(g_{\ti{\phi}_1 \,xv_1} + \cdots
  g_{\ti{\phi}_k\, xv_k}\right)(U(\g)).
\end{eqnarray*} 
By induction we find $\ti{\phi}_1=\cdots=\ti{\phi}_k=0$. Inserting in equation
(\ref{mbeq}), with (\ref{fmu}) follows 
\begin{eqnarray*}
  \{0\}= \left(g_{\phi_{k+1} \ti{v}_{k+1}}+ \cdots + g_{\phi_n \ti{v}_n}\right)(U(\g)).
\end{eqnarray*}
By induction
$\phi_{k+1}=\cdots=\phi_{n}=0$. Since $\ti{\phi}_1=\cdots=\ti{\phi}_k=0$ this also implies
$\phi_1=\cdots=\phi_k=0$.

Now assume that for every $x\in U(\g)$ either
 $x v_1,\,\ldots,\, xv_n$ are $S(V)$-linearly independent, or all
 these elements are zero. 
 Then for every $x\in U(\g)$ we have $xv_1=0$ if and only if $xv_2=0$. Therefore we
 get a well defined $\g$-equivariant non-zero linear map $\al: U(\g)v_1\to
 U(\g)v_2 $ by $\al(xv_1):=xv_2$. Since $V$ is irreducible
 $U(\g)v_1=U(\g)v_2=V$. We have found an element $\al\in S(V)$ such
 that $\al(v_1)=v_2$, contradicting the $S(V)$-linear independence of $v_1$, $v_2$.
\qed
\subsection{Jordan-Chevalley decompositions for elements of Lie(M) and $M$}
Recall the additive Jordan-Chevalley decomposition for locally finite endomorphisms stated in Section \ref{Pre}.
We call an element $x\in Lie(M)$ {\it locally finite}, resp. {\it semisimple},
  resp. {\it locally nilpotent}, if it acts locally finite, resp. semisimple,
  resp. locally nilpotent on every module contained in $\mathcal C$. 
\begin{Theorem}\label{JChadd} Let $x\in Lie(M)$ be locally finite. Suppose that for every module $V$ contained in $\mathcal C$ the semisimple endomorphism 
$s_V$ and the locally nilpotent endomorphism $n_V$ of the Jordan-Chevalley decomposition of $x_V$ are already in $End_{V^{du}}(V)$. Set
\begin{eqnarray*} 
  s:=(s_V)_{V\;an\; object\; of\; {\mathcal C}} \quad \mb{ and }\quad n:=(n_V)_{V\;an\; object\; of\; {\mathcal C}}. 
\end{eqnarray*}
Then $s$ is a semisimple element of $Lie(M)$ and $n$ is a locally nilpotent
element of $Lie(M)$. They satisfy $x=s+n$ and $[s,n]=0$, and are uniquely determined by these properties.
\end{Theorem}
\begin{Remark} The condition of the theorem on the locally finite element $x\in Lie(M)$ is satisfied if $V^{du}=V^*$ for all modules $V$ contained in $\mathcal C$. 
It is possible to show that the condition is also satisfied if $x_V^{du}$ acts locally finite on $V^{du}$ for all modules $V$ contained in 
$\mathcal C$. (Reduce to an algebraic closed field. Then use the generalized eigenspaces of $x_V$ and $x^{du}_V$ to show that $s^{du}_V$ is the semisimple part of the 
Jordan-Chevalley decomposition of $x^{du}_V$.)
\end{Remark}

\Proof (a) Let $U$, $\ti{U}$ be $\F$-linear spaces. Let $x\in End(U)$, $\ti{x}\in End(\ti{U})$ locally finite endomorphisms with Jordan-Chevalley
decompositions $x=s+n$, $\ti{x}=\ti{s}+\ti{u}$. Let $\al: U\to\ti{U}$ be a
linear map such that $\al\circ x = \ti{x}\circ \al$. We show
\begin{eqnarray*}
   \al\circ s = \ti{s}\circ \al \quad\mb{ and }\quad \al\circ n = \ti{n}\circ \al.
\end{eqnarray*}
It is sufficient to show the first equation. Since $x$ is locally finite it is sufficient to assume $U$
finite dimensional. 
Because of $\al\circ x = \ti{x}\circ \al$, the image $\al(u)$ is $\ti{x}$-invariant. It is sufficient to assume $\al:U\to\ti{U}$ surjective.
By Proposition 19 of Section 7.1 in \cite{MoPi} it is sufficient to assume $\F$ to be algebraic closed. 
Let $U=\bigoplus_{\la\in\F}U^\la$ and $\ti{U}=\bigoplus_{\la\in\F}\ti{U}^\la$ be the generalized eigenspace
decompositions of $x$ and $\ti{x}$. With $\al\circ x = \ti{x}\circ \al$ it follows 
$\al(U^\la)=\ti{U}^\la$ for all $\la\in\F$. 
By Proposition 20 of Section 7.1 in \cite{MoPi} the endomorphism $s$ acts on $U^\la$ diagonally by the eigenvalue $\la$. Similarly, 
$\ti{s}$ acts on $\ti{U}^\la$ diagonally by the eigenvalue $\la$. Therefore $\al\circ s = \ti{s}\circ \al$.\vspace*{1ex} 

(b) Now let $x \in Lie(M)$ be locally finite. Immediately form (a) and from
the condition of the theorem follows $s,n\in Nat$.
Recall that for an element $y\in Nat$ it is equivalent:
\begin{itemize}
\item[(i)] $y\in Lie(M)$.
\item[(ii)] There exists a derivation $\pi(1\otimes y):\FK{M}\to \FK{M}$
such that 
\begin{eqnarray*}
  \pi(1\otimes y)f_{\phi v}= f_{\phi\, y_V v} \quad\mb{ for all }\quad \phi\in V^{du},\;v\in V,\; V\mb{ an object of }{\mathcal C}.
\end{eqnarray*}
\end{itemize}

Since $x$ is locally finite also $\pi(1\otimes x)$ is locally finite. Let $\pi(1\otimes s)$ be the semisimple part, and $\pi(1\otimes n)$ the locally nilpotent part of the
Jordan-Chevalley decomposition of $\pi(1\otimes x)$. By Proposition 26 of Section 7.1 in \cite{MoPi} the endomorphisms $\pi(1\otimes s)$ and $\pi(1\otimes n)$ are derivations
of $\FK{M}$. 
Now let $V$ be a $\g$-module contained in $\mathcal C$, and $\phi\in V^{du}$. The linear map
\begin{eqnarray*}
  \Phi:\, V \mapsto \FK{M}\quad \mb{ defined by }\quad \Phi(v) :=f_{\phi v} 
\end{eqnarray*}
satisfies $\pi(1\otimes x)\circ\Phi=\Phi\circ x_V$. By (a) we get $\pi(1\otimes s)\circ\Phi=\Phi\circ s_V $ and $\pi(1\otimes n)\circ\Phi=\Phi\circ n_V $,
from which follows
\begin{eqnarray*}
   \pi(1\otimes s) f_{\phi v}= f_{\phi\, s_V v} \quad\mb{ and }\quad \pi(1\otimes n) f_{\phi u}= f_{\phi\, n_V v}\quad\mb{ for }\quad v\in V.
\end{eqnarray*}
Therefore $s,n\in Lie(M)$.
\qed 

Let $e\in M$ be an idempotent. Then
\begin{eqnarray*}
   M_e:=e M e=\Mklz{m\in M}{ me=em =m }
\end{eqnarray*}
is a closed submonoid of $M$ with unit $e$. In particular $M_1=M$. 
Denote by $M_e^\times$ the unit group of the monoid $M_e$. It is easy to check that for $e\neq e'$ we have 
$M_e^\times\cap M_{e'}^\times = \emptyset$. In the next theorem we give a
multiplicative Jordan-Chevalley decomposition for certain elements of 
\begin{eqnarray*}
   \dot{\bigcup_{e\;an\;idempotent \;of\; M}} M_e^\times.
\end{eqnarray*}

Recall the multiplicative Jordan-Chevalley decomposition of Section \ref{Pre}.
We call an element $x\in M$ {\it locally finite}, resp. {\it semisimple},
  resp. {\it locally weak unipotent}, if it acts locally finite, resp. semisimple,
  resp. locally weak unipotent on every module contained in $\mathcal C$. 
\begin{Theorem}\label{JChmult}Let $e\in M$ be an idempotent, and $m\in M_e^\times$ locally
  finite. Suppose that for every module $V$ contained in $\mathcal C$ the semisimple endomorphism 
$s_V$ and the locally weak unipotent endomorphism $u_V$ of the
Jordan-Chevalley decomposition of $m_V\in ((End_{V^{du}}(V))_{e_V})^\times$ are already in $End_{V^{du}}(V)$. Set
\begin{eqnarray*} 
  s:=(s_V)_{V\;an\; object\; of\; {\mathcal C}}\quad\mb{ and }\quad
  u:=(u_V)_{V\;an\; object\; of\; {\mathcal C}}. 
\end{eqnarray*}
Then $s$ is a semisimple element of $M_e^\times$ and $u$ is a locally weak
unipotent element of $M_e^\times$. They satisfy $m=su=us$ and are uniquely determined by these properties.
\end{Theorem}
\begin{Remark} The condition of the theorem on the locally finite element $m\in M_e$ is satisfied if $V^{du}=V^*$ for all modules $V$ 
contained in $\mathcal C$. It is possible to show that the condition is also satisfied if
$m_V^{du}$ acts locally finite on $V^{du}$ for all modules $V$ contained in $\mathcal C$.
\end{Remark}

\Proof The essential part of the proof of $s,u\in Nat$ is a variant of the corresponding
proof for the additive Jordan-Chevalley decomposition of $Lie(M)$.

From $m_{V_0}=id_{V_0}$ follows $s_{V_0}=id_{V_0}$ and $u_{V_0}=id_{V_0}$. Now
let $V$, $W$ be $\g$-modules contained in $\mathcal C$. We have $m_V\otimes
m_W=m_{V\otimes W}$. To prove the analogous equations for $s$ and $n$ it is
sufficient to show that $s_V\otimes s_W$ is the semisimple part, and $u_V\otimes
u_W$ is the locally weak unipotent part of the Jordan-Chevalley decomposition of $m_V\otimes
m_W\in \left((End_{(V\otimes W)^{du}}(V\otimes W))_{e_{V\otimes W}}\right)^\times$.
This follows easily by the results stated in Section \ref{Pre}.
\qed
\subsection{The case where the $\g$-modules contained in $\mathcal C$ are finite dimensional}
Let $M$ be the Tannaka monoid with coordinate ring $\FK{M}$ and Lie algebra
$Lie(M)$ obtained from a pair of categories $\mathcal C$, ${\mathcal C}^{du}$.
In general the multiplication map of the Tannaka monoid $M$ does not induce a comultiplication of
its coordinate ring of matrix coefficients $\FK{M}$. 
Assigning to an element of $M$ its evaluation homomorphism, $M$ embeds as a
set into $\Spm \FK{M}$. In general this map is not surjective. Its image is
described in part (a) of the following lemma. Similarly, assigning to an element of $Lie(M)$ its corresponding derivation of $\FK{M}$ in $1\in M$, 
the linear space $Lie(M)$ embeds into $Der_1(\FK{M})$. In general, also this map is not surjective. Its image is described in part (b) of the following Lemma. 

\begin{Lemma}\label{hd} Let $M$ be the Tannaka monoid obtained 
from the pair of categories $\mathcal C$, ${\mathcal C}^{du}$. Let $\FK{M}$ be its coordinate ring of matrix coefficients. 
\begin{itemize}
\item[(a)] Let $\al\in\Spm\FK{M}$. Then $\al$ is the evaluation homomorphism of
an element of $M$ if and only if for every object $V$ contained in $\mathcal
C$ there exists an endomorphism $\al_V\in End_{V^{du}}(V)$ such that
\begin{eqnarray}\label{alev}
   \phi(\al_V v)=\al(f_{\phi v})\quad \mb{ for all }\quad v,\in V,\;\phi\in V^{du}.
\end{eqnarray}
\item[(b)] Let $\delta\in Der_1(\FK{M})$. Then $\delta$ is the derivation of an element of $Lie(M)$ if and only if for every object $V$ 
contained in $\mathcal C$ there exists an endomorphism $\delta_V\in End_{V^{du}}(V)$ such that
\begin{eqnarray*}
   \phi(\delta_V v)=\delta(f_{\phi v})\quad \mb{ for all }\quad v,\in V,\;\phi\in V^{du}.
\end{eqnarray*}
\end{itemize}
\end{Lemma}

\Proof We only show part (a) of the proposition, part (b) is
proved similarly. We only prove the 'if' part of this statement, the 'only if' part is trivial. 
Let $V$, $W$ be $\g$-modules contained in $\mathcal C$. Let $\eta:V\to
W$ be a homomorphism of $\g$-modules. Then for all $v\in V$ and $\psi\in W^*$
we find
\begin{eqnarray*}
   \psi(\eta\al_V v) =  (\eta^*\psi)(\al_V v) = 
                       \al(f_{\eta^*\psi\,v})= \al(f_{\psi\,\eta v})= \psi(\al_W\eta v).
\end{eqnarray*}
It follows $\eta\circ\al_V=\al_W\circ\eta$.

Let $V$, $W$ be $\g$-modules contained in $\mathcal C$. For all $v\in V$, $w\in W$, $\phi\in V^*$, $\psi\in W^*$
we have 
\begin{eqnarray*}
  (\phi\otimes\psi)\left(\al_{V\otimes W}(v\otimes
  w)\right)= \al(f_{\phi\otimes\psi\,v\otimes w})=
           \al(f_{\phi v}) \al(f_{\psi w})=\\
   \phi(\al_V v)\psi(\al_W w)= (\phi\otimes\psi)\left((\al_v\otimes\al_w)(v\otimes w)\right).
\end{eqnarray*}
Therefore $\al_{V\otimes W}=\al_v\otimes\al_w$. Let $V_0$ be a one-dimensional trivial $\g$-module. If $v_0\in V_0$ and $\phi_0\in V_0^{du}$ then 
$\al(f_{\phi_0 v_0})=\al(\phi_0(v_0)1)=\phi_0(v_0)$. It follows $\al_{V_0}=id_{V_0}$.

We have shown that $(\al_V)_{V\;an\;object\;of\;{\mathcal C}}$ is an element of
$M$. By (\ref{alev}) its evaluation homomorphism coincides with $\al$.
\End

If the $\g$-modules contained in $\mathcal C$ are finite dimensional then the Tannaka reconstruction of this paper coincides with the classical 
Tannaka reconstruction as described in \cite{DeMi}. The following two theorems can be obtained by specializing the results of \cite{DeMi}. 
It is also easy to give direct proofs, which we include here in the ArXiv-version. 
Note that if the $\g$-modules contained in $\mathcal C$ are finite dimensional then ${\mathcal C }^{full}$ is the only category of duals of 
$\mathcal C$. 

\begin{Theorem}\label{oneparameterfinite2} Suppose that every $\g$-module contained in $\mathcal C$ is finite dimensional. Let $M$ be the Tannaka 
monoid obtained from $\mathcal C$ and ${\mathcal C}^{full}$, let $\FK{M}$ be its
coordinate ring of matrix coefficients, and let $Lie(M)$ be its Lie algebra. 
Then the monoid structure of $M$ induces in the natural way a bialgebra structure on the coordinate ring $\FK{M}$.
As a set, $M$ identifies with $\Spm\FK{M}$. As a linear space, $Lie(M)$ identifies with $Der_1(\FK{M})$.
\end{Theorem}

\Proof (a) The comorphism of the multiplication map of $M$ exists: Let $V$ be a $\g$-module contained in $\mathcal C$. Let $v\in V$
and $\phi\in V^*$. Choose a base $v_1,\,\ldots,\,v_n$ of $V$, and let 
$\phi_1,\,\ldots,\,\phi_n$ be the dual base of $V^*$. Then for $m,\ti{m}\in M$ we get
\begin{eqnarray*}
  f_{\phi v}(m\ti{m})=\sum_{i=1}^n f_{\phi v_i}(m)f_{\phi_i v}(\ti{m}).
\end{eqnarray*}

The algebra $\FK{M}$ gets the structure of a bialgebra by taking as comultiplication
the comorphism of the multiplication map of $M$, and as counit the evaluation map in the unit 
of $M$. \vspace*{1ex}

(b) Let $\al\in\Spm\FK{M}$. Let $V$ be a 
$\g$-module contained in $\mathcal C$. Choose a base $v_1,\,\ldots,\,v_n$ of $V$, and let
$\phi_1,\,\ldots,\,\phi_n$ be its dual base. Define an endomorphism $\al_V \in
End(V)$ by
\begin{eqnarray*}
  \al_V v:=\sum_{i=1}^n \al(f_{\phi_i\, v}) v_i.
\end{eqnarray*} 
Then for all $v\in V$ and $\phi\in V^*$ we have
\begin{eqnarray*}
   \phi(\al_V v)= \sum_{i=1}^n \al(f_{\phi_i\, v}) \phi(v_i)
                       = \al(f_{\sum_i \phi_i \phi(v_i)\,v })
                       = \al(f_{\phi\, v}).
\end{eqnarray*}
By part (a) of Lemma \ref{hd} the homomorphism $\al$ is the evaluation homomorphism of an element of $M$.\vspace*{1ex}

(c) Similarly, with part (b) of Lemma \ref{hd} follows that $Lie(M)$ identifies with $Der_1(\FK{M})$. 
\qed

\begin{Theorem}\label{oneparameterfinite} Suppose that every $\g$-module contained in $\mathcal C$ is finite dimensional. Suppose that in 
addition for every $\g$-module $V$ contained in $\mathcal C$ also its dual $\g$-module $V^*$ is contained in $\mathcal C$. 
Then the Tannaka monoid $M$ obtained from $\mathcal C$ and ${\mathcal C}^{full}$ is a group. 
The group structure of $M$ induces in the natural way a Hopf algebra structure
on the coordinate ring of matrix coefficients $\FK{M}$.
\end{Theorem} 

\Proof For an endomorphism $\beta\in End(V^*)$ of a finite dimensional linear space $V$ its adjoint $\beta^*\in End(V)$ exists. 
For an element $m\in M$ define 
\begin{eqnarray*}
   (m^{inv})_V:= (m_{V^*})^* \quad \mb{ where } \quad V \mb{ is an object of }{\mathcal C}.
\end{eqnarray*} 
For every homomorphism
$\al:V\to W$ of $\g$-modules, $V$, $W$ objects of $\mathcal C$, we have
\begin{eqnarray*}
  \al\circ  (m_{V^*})^* =(m_{V^*}\circ \al^*)^* = (\al^*\circ m_{W^*})^* =
  (m_{W^*})^*\circ \al.
\end{eqnarray*}
Furthermore, for all $\g$-modules $V$, $W$ contained in $\mathcal C$ we have
\begin{eqnarray*}
 (m^{inv})_{V\otimes W}=(m_{(V\otimes W)^*})^*=(m_{V^*}\otimes m_{W^*})^*=(m^{inv})_V\otimes(m^{inv})_W. 
\end{eqnarray*}
If $V_0$ is a trivial one-dimensional $\g$-module, then also its dual $V_0^*$ is a trivial one-dimensional $\g$-module. 
It follows $m_{V_0}^{inv}=(m_{V_0^*})^*=(id_{V_0^*})^*=id_{V_0}$.
Therefore we get a map $\mb{}^{inv}:M\to M$. The comorphism exists, because for every $\g$-module $V$ contained in $\mathcal C$, $v\in V$, and
$\phi\in V^*$ we have
\begin{eqnarray*}
  f_{\phi v}(m^{inv})=\phi((m_{V^*})^*v)=(m_{V^*}\phi)(v)=f_{v\phi}(m)
  \quad\mb{ for all } \quad m\in M.
\end{eqnarray*}
 
Let $m\in M$. Since the map $\delta:V^*\otimes V\to\F$ defined by $\delta(\phi\otimes v):=\phi(v)$, where $\phi\in V^*$, $v\in V$ is 
$\g$-equivariant it is also $M$-equivariant. Therefore, we have
\begin{eqnarray*}
   \phi(v)= m\left(\phi(v)\right)=(m\phi)(mv)=\phi(m^{inv}mv)\quad\mb{ for all }\quad \phi\in V^*,\;v\in V. 
\end{eqnarray*}
It follows $m^{inv}m=1$. Let $v_1,\,\ldots,\,v_n$ a base of $V$, and let
$\phi_1,\,\ldots,\,\phi_n$ the dual base of $V^*$. Since the map $\mu:\F\to
V^*\otimes V$ defined by $\mu(c):=c\sum_{i=1}^n\phi_i\otimes v_i$, where
$c\in\F$ is $\g$-equivariant it is also $M$-equivariant. Therefore, we have
\begin{eqnarray*}
  \sum_{i=1}^n \phi_i\otimes v_i= \mu(1)=\mu(m 1)=m\mu(1)=\sum_{i=1}^n m \phi_i\otimes m v_i.
\end{eqnarray*} 
From this equation we get
\begin{eqnarray*}
\psi(v)=\sum_{i=1}^n\phi_i(v)\psi(v_i)=\sum_{i=1}^n(m\phi_i)(v)\psi(m v_i)=
\sum_{i=1}^n \phi_i (m^{inv}v)\psi(m v_i) =\psi(m m^{inv}v)
\end{eqnarray*}
for all $\psi\in V^*$, $v\in V$. It follows $m m^{inv}=1$.

The coordinate ring $\FK{M}$ equipped with the bialgebra structure of $\FK{M}$
together with the comorphism of the map $inv$ as antipode is a Hopf algebra. 
\qed

Actually, we are not only interested to obtain the Tannaka monoid, its coordinate ring, and its Lie algebra. We also want to realize the Lie 
algebra from which we started as a subalgebra of the Lie algebra of the Tannaka monoid.

\begin{Theorem}\label{fdgood} Suppose that every $\g$-module contained in $\mathcal C$ is finite-dimensional. Then the pair 
$\mathcal C$, ${\mathcal C}^{full}$ is very good for integrating $\g$.
\end{Theorem}

\Proof (a) We first show that the pair $\mathcal C$, ${\mathcal C}^{full}$ is
good for integrating $\g$. Recall the algebra of matrix coefficients $\FK{U(\g)}$ on $U(\g)$ introduced in Theorem \ref{UnivEin} (a). 
Recall that $\Spm \FK{U(\g)}$ denotes the set of homomorphisms of algebras from $\FK{U(\g)}$ to $\F$. In particular, the evaluation map 
$\al_1:\FK{U(\g)}\to \F$ in $1\in U(\g)$ is an element of $\Spm \FK{U(\g)}$. Every element $g\in \FK{U(\g)}$ gives a function on $\Spm
\FK{U(\g)}$ by $g(\al):=\al(g)$, $\al\in \Spm \FK{U(\g)}$. Identify
$\FK{U(\g)}$ with the corresponding coordinate ring on $\Spm \FK{U(\g)}$. 

Let $\al\in \Spm \FK{U(\g)}$ and $V$ be a $\g$-module contained in $\mathcal C$. Choose a base $v_1,\,\ldots,\,v_n$ of $V$. Let 
$\phi_1,\,\ldots,\,\phi_n$ be the dual base of $V^*$. Define a map $\al_V:V\to V$ by
\begin{eqnarray*}
  \al_V v:=\sum_{i=1}^n \al(g_{\phi_i v}) v_i,\quad\mb{ where }\quad v\in V.
\end{eqnarray*}
The map $\al_V$ is independent of the choice of the base of $V$ and linear, since
\begin{eqnarray*}
   \phi(\al_V v)= \sum_{i=1}^n \al(g_{\phi_i v}) \phi(v_i)=\al(g_{\sum_{i=1}^n
   \phi(v_i)\phi_i \,v}) =\al(g_{\phi v})
\end{eqnarray*}
for all $\phi\in V^{du}$. 
Fix an element $\al\in \Spm \FK{U(\g)}$. Let $V$, $W$ be $\g$-modules contained in $\mathcal C$,
and let $\gamma:V\to W$ a $\g$-equivariant linear map. Then for all $v\in V$ and $\psi\in W^*$ we find
\begin{eqnarray*}
   \psi(\al_W\gamma v)=\al(g_{\psi\,\gamma(v)})=\al(g_{\gamma^{du}(\psi)\,v})=\gamma^{du}(\psi)(\al_V v) = \psi( \gamma \al_V v).
\end{eqnarray*}
It follows $\al_W\circ\gamma=\gamma\circ \al_V$. 
If $V_0$ is a one-dimensional trivial $\g$-module then trivially $\al_{V_0}=id_{V_0}$. Now let $V$, $W$ be $\g$-modules contained 
in $\mathcal C$. For all  $v\in V$, $w\in W$, and for all $\phi\in V^{du}$,
$\psi\in W^{du}$ we find by using equation (\ref{Umult}) of the proof of
Theorem \ref{UnivEin}:
\begin{eqnarray*}
     (\phi\otimes\psi)\left(\al_{V\otimes W}(v\otimes w)\right)=
     \al(g_{\phi\otimes\psi\,v\otimes w})=\al(g_{\phi v}\cdot g_{\psi w})=\al(g_{\phi v})\al(g_{\psi w})\qquad\qquad\\
    =\phi(\al_V v)\psi(\al_W w)=(\phi\otimes\psi)\left( (\al_V\otimes \al_W)(v\otimes w)\right).
\end{eqnarray*}
It follows $\al_{V\otimes W}=\al_V\otimes \al_W$.
Therefore we get a map $i:\Spm \FK{U(\g)}\to M$ defined by $i(\al):=(\al_V)_{V\;an\;object\;of\;{\mathcal C}}$, $\al\in \Spm \FK{U(\g)}$.

For all $\al\in \Spm \FK{U(\g)}$ and all $\g$-modules $V$ contained in $\mathcal C$, 
$v\in V$, and $\phi\in V^*$ we have $f_{\phi v}(i(\al))=\phi(\al_V v)=g_{\phi v}(\al)$. This shows that $i$ is
injective. Via the map $i$ the coordinate ring of $M$ restricted to $i(\Spm \FK{U(\g)})$ identifies with $\FK{U(\g)}$.

The evaluation map $\ti{\delta}_x:\FK{U(\g)}\to \F$ in an element $x\in\g\subseteq
U(\g)$ is a derivation of $\FK{U(\g)}$ in $\al_1\in \Spm \FK{U(\g)}$. 
Let $\delta:\FK{M}\to\F$ be the concatenation of the maps
\begin{eqnarray*}
   \FK{M}\to\FK{i(D)}\cong \FK{U(\g)}\to \F,
\end{eqnarray*} 
the first the restriction map onto $i(D)$, the second the comorphism $i^*$, and the last the map $\ti{\delta}_x$. Then $\delta$
is a derivation in $1\in M$ such that $\delta(f_{\phi v})=\phi(xv)$ for all
$\phi\in V^*$, $v\in V$, and $V$ a $\g$-module contained in $\mathcal
C$. Therefore $x\in Lie(M)$.\vspace*{1ex} 

 (b) To show that the pair $\mathcal C$, ${\mathcal C}^{full}$ is very good for integrating $\g$ it is sufficient to show that the unit group 
$M^\times$ is dense in $M$, compare Proposition \ref{Defrel}.
Let $\widetilde{\mathcal C}$ be the full, wide subcategory of the category of
$\g$-modules, closed under taking finite sums, 
tensor products, submodules, which is  
generated by the $\g$-modules contained in $\mathcal C$ and the dual $\g$-modules.
By Theorem \ref{oneparameterfinite} the Tannaka monoid associated to the pair $\widetilde{\mathcal C}$, 
${\widetilde{\mathcal C}}^{full}$ is a group, which we denote by
$G$. Obviously, we get a homomorphism of groups by: 
\begin{eqnarray*}
  i:\,G&\to& M^\times\\
   (g_V)_{V\;an\;object\;of\;\widetilde{\mathcal C}}&\to&  (g_V)_{V\;an\;object\;of\;{\mathcal C}}
\end{eqnarray*}
Let $V$ be an object of $\mathcal C$. Let $v\in V$, and $\phi\in V^*$ such that $\phi(M^\times v)=\{0\}$. Then also $\phi(G v)=\{0\}$. By part (a), the
pair $\widetilde{\mathcal C}$, ${\widetilde{\mathcal C}}^{full}$ is good for
integrating $\g$. It follows $\phi(U(\g))v)=\{0\}$ by Theorem \ref{UnivEin}.  
By part (a) also the
pair $\mathcal C$, ${\mathcal C}^{full}$ is good for
integrating $\g$. Applying Theorem \ref{UnivEin} once more we get $\phi(M v)=\{0\}$.
\qed
\subsection{An embedding theorem}

As before let $M$ be the Tannaka monoid associated to a category of $\g$-modules $\mathcal C$ and a category of duals ${\mathcal C}^{du}$. 
Let $\FK{M}$ be its coordinate ring and $Lie(M)$ its Lie algebra.
Now suppose that the Lie subalgebra $\bf s$ of $\g$ acts locally finite on every $\g$-module contained in $\mathcal C$. 
Define ${\mathcal C}({\bf s})$ to be the category of $\bf s$-modules, whose objects are the $\bf s$-modules isomorphic to some finite-dimensional 
$\bf s$-submodule of some $\g$-module contained in $\mathcal C$. (As morphisms between two $\bf s$-modules we take all $\bf s$-equivariant maps.)
The Tannaka reconstruction may applied to the pair ${\mathcal C}({\bf s})$ and ${\mathcal  C}({\bf s})^{full}$, 
giving the Tannaka monoid $M({\bf s})$ with coordinate ring $\FK{M({\bf s})}$, and the Lie algebra $Lie(M({\bf s}))$.

Every object $V$ of $\mathcal C$ is a locally finite $\bf s$-module, whose finite-dimensional $\bf s$-submodules are contained in 
${\mathcal C}({\bf s})$. The monoid $M({\bf s})$ as well as its Lie algebra $Lie(M({\bf s}))$ act in the obvious way on $V$. For $m\in M({\bf s})$
resp. $x\in Lie(M({\bf s}))$ we denote by $m_V\in End(V)$ resp. $x_V\in End(V)$ the corresponding endomorphism of $V$.

\begin{Theorem}\label{embedding} Let $\bf s$ be a Lie subalgebra of $\g$ which acts locally finite on every $\g$-module contained in $\mathcal C$. 
With the notations introduced before set
\begin{eqnarray*}
     D    &:=& \Mklz{m\in M({\bf s})}{m_V \in End_{V^{du}}(V) \mb{ for all objects } V \mb{ of }{\mathcal C}},\\
  {\bf d} &:=& \Mklz{x\in Lie(M({\bf s}))}{x_V \in End_{V^{du}}(V) \mb{ for all objects } V \mb{ of }{\mathcal C}}.
\end{eqnarray*} 
Then $D$ is a submonoid of $M({\bf s})$, and ${\bf d}$ is a Lie subalgebra of
$Lie(M({\bf s}))$. Both act locally finite on the objects of $\mathcal C$. Furthermore:
\begin{itemize}
\item[(a)] The map
\begin{eqnarray*}
  i:\; D &\to& M\\
   d  &\to& (d_V)_{V\;object\;of\;{\mathcal C}}
\end{eqnarray*}
is a morphism of monoids which is also a closed embedding of sets with coordinate rings.
\item[(b)] Suppose that $D$ is dense in $M({\bf s})$. Then the map
\begin{eqnarray*}
 Lie(i):\;  {\bf d} &\to & Lie(M)\\
    x  &\to& (x_V)_{V\;object\;of\;{\mathcal C}}
\end{eqnarray*}
is an embedding of Lie algebras with image $Lie(i(D))$. In particular, we have ${\bf  s}\subseteq Lie(i(D))\subseteq Lie(M)$. 
\end{itemize}
\end{Theorem}

\begin{Remark} If $V^{du}=V^*$ for all objects of $\mathcal C$, then $D=M({\bf s})$ and ${\bf d}= Lie(M({\bf s}))$.
\end{Remark}

\Proof Obviously $D$ is a submonoid of $M({\bf s})$ and $\bf d$ a Lie
subalgebra of $Lie(M({\bf s}))$. To (a): Every $\g$-module contained in $\mathcal C$ restricts to a locally
finite $\bf s$-module, whose finite dimensional $\bf s$-submodules are contained in ${\mathcal C}({\bf s})$. 
Every morphism of $\g$-modules contained in $\mathcal C$ restricts to a morphism of locally finite $\bf s$-modules. 
Therefore $M({\bf s})$ acts on the $\g$-modules contained in $\mathcal C$,
compatible with morphisms and tensor products. Since $M({\bf s})$ acts already
faithfully on the finite dimensional ${\bf s}$-submodules of the $\g$-modules contained in
$\mathcal C$, it acts faithfully on the $\g$-modules contained in
$\mathcal C$. This shows that the map $i:D\to M$ is a well defined injective
morphism of monoids.

If $v\in V$, $\phi\in V^{du}$, and $V$ be an object of $\mathcal C$ then
\begin{eqnarray*}
  \phi(i(d)v)=\phi\res{U({\bf s})v}(d v) \quad\mb{ for all }d\in D.
\end{eqnarray*}
Here $U({\bf s})v$ is a $\bf s$-module contained in ${\mathcal C}({\bf  s})$. Since $V^{du}$ is point separating on $V$ it restricts to the full 
dual on the finite-dimensional subspace $U({\bf s})v$ of $V$. Therefore the
comorphism $i^*$ exists and is surjective. 

Denote the restriction of the map $i:D\to M$ onto its image also by $i$. Then $i:D\to i(D)$ is an isomorphism of sets with coordinate ring. Now let
$\ti{m}\in \overline{i(D)}$ and let $\al_{\ti{m}}:\FK{ \overline{i(D)}}\to\F$ the corresponding evaluation homomorphism. Consider the 
concatenation of the homomorphisms of algebras
\begin{eqnarray*}
  \FK{M({\bf s})}\to\FK{D}\cong \FK{i(D)}\cong\FK{\overline{i(D)}}\to\F,
\end{eqnarray*} 
the first the restriction map, the second given by the inverse of the
isomorphism $i^*:\FK{i(D)}\to \FK{D}$, the third given by the inverse of the
restriction map, which exists since $i(D)$ is dense in $\overline{i(D)}$, and
the last the evaluation morphism $\al_{\ti{m}}$. By Theorem \ref{oneparameterfinite2} this
homomorphism is an evaluation homomorphism in an element $m\in M({\bf s})$. Therefore
\begin{eqnarray*}
  \phi(\ti{m}_V v)=\phi\res{U({\bf s})v}(m_V v)\quad\mb{ for all }\quad \phi\in V^{du},\;v\in V,\;V\mb{ an object of }{\mathcal C}.
\end{eqnarray*}
It follows $\ti{m}_V=m_V$ for all objects of $\mathcal C$. Therefore $m\in D$ and $\ti{m}=i(m)$.\vspace*{1ex}

To (b): Let $x\in{\bf d}$. let $\ti{\delta}_x$ the corresponding derivation in
$1\in M({\bf s})$ of $\FK{M({\bf s})}$. Consider the concatenation $\delta$ of the maps
\begin{eqnarray*}
  \FK{M}\to \FK{i(D)}\cong \FK{D}\cong\FK{M({\bf s})}\to\F,
\end{eqnarray*}
the first the restriction map, the second the comorphism $i^*$, the third the
inverse of the restriction map, which exists since $D$ is dense in $M({\bf
  s})$, and the last the derivation $\ti{\delta}_x$. It is easy to check that
$\delta$ is a derivation in $1\in M$ and 
\begin{eqnarray*}
  \delta(f_{\phi v})=\ti{\delta}(f_{\phi\res{U({\bf s})v} v}) =\phi(x_V v)\quad\mb{ for all }\quad \phi\in
  V^{du},\;v\in V,\;V \mb{ an object of }{\mathcal C}.
\end{eqnarray*}
Therefore $\delta=\delta_{i(x)}$. Obviously the kernel of $\delta$ contains
the vanishing ideal $I(i(D))$. It follows $i(x)\in Lie(i(D))$.

By Theorem \ref{fdgood} the Lie algebra $\bf s$ identifies with a subalgebra
of $Lie(M({\bf s}))$. The action of an element $x\in {\bf s}$ as an element of
$Lie(M({\bf s}))$ on a $\g$-module $V$ contained in $\mathcal C$ coincides with the
action of $x$ as an element of $\g$ on $V$. It follows that $\bf s$ identifies
with a subalgebra of $\bf d$. This subalgebra is mapped by $Lie(i)$ onto the
Lie subalgebra of $Nat$, given by $\bf s$.

Now let $y\in Lie(i(D))$. Let $\ti{\delta}_y$ the corresponding derivation in
$1\in i(D)$ of $\FK{i(D)}$. Consider the concatenation of the maps
\begin{eqnarray*}
       \FK{M({\bf s})}\cong \FK{D}\cong\FK{i(D)}\to\F,
\end{eqnarray*}
the first the restriction map, the second the inverse of the comorphism $i^*$,
the third the derivation $\ti{\delta}_y$. This concatenation is a derivation
of $\FK{M({\bf s})}$ in $1\in M({\bf s})$. By Theorem \ref{oneparameterfinite2} it is a 
derivation corresponding to an element $x\in Lie(M({\bf s}))$. Therefore
\begin{eqnarray*}
  \phi(y_V v)=\phi\res{U({\bf s})v}(x_V v)\quad\mb{ for all }\quad \phi\in V^{du},\;v\in V,\;V\mb{ an object of }{\mathcal C}.
\end{eqnarray*}
It follows $y_V=x_V$ for all objects of $\mathcal C$. Therefore $x\in {\bf d}$ and $y=\left(Lie(i)\right)(x)$. 
\qed

\begin{Remark}\label{construct} To apply the last theorem in concrete situations, we have to determine the monoid $M({\bf s})$, its coordinate 
ring $\FK{M({\bf s})}$, and its Lie algebra $Lie(M({\bf s}))$, as well as the actions of $M({\bf s})$ and $Lie(M({\bf s}))$ on the objects of 
$\mathcal C$. There is a direct possibility which is sometimes useful. 

The coordinate ring $\FK{M({\bf s})}$ is isomorphic to the the algebra $\FK{U({\bf s})}$ of matrix coefficients on $U({\bf s})$. 
For every $\g$-module $V$ contained in $\mathcal C$, the dual $V^{du}$ is point separating on $V$. It restricts to the full 
dual on every finite-dimensional $\bf s$-subspace of $V$. Therefore, the restriction of the algebra of matrix coefficients 
$\FK{U(\g)}$ on $U(\g)$ onto $U({\bf s})\subseteq U(\g)$ coincides with $\FK{U({\bf s})}$. 

By Theorem \ref{oneparameterfinite2} we may identify the set $M({\bf s})$ with $\Spm \FK{U({\bf s})}$. 
If $\al\in \Spm \FK{U({\bf s})}$, and $V$ is an object of $\mathcal C$, the corresponding endomorphism $\al_V\in End(V)$ is determined by
\begin{eqnarray*}
  \phi(\al_V v)=\al(g_{\phi v} \res{U({\bf s})} ) \quad\mb{ for all }\quad \phi\in V^{du},\;v\in V. 
\end{eqnarray*}
The monoid structure on $\Spm \FK{U({\bf s})}$ is the uniquely determined monoid structure, such that the maps  
$\Spm \FK{U({\bf s})}\to End(V)$, $V$ an object of $\mathcal C$, are morphisms of monoids.

Similarly, by Theorem \ref{oneparameterfinite2} we may identify the linear space $Lie(M({\bf s}))$ with $Der_1(\FK{U({\bf s})})$.
(For this notation note, that the evaluation homomorphism of $\FK{M({\bf s})}$ in $1\in M({\bf s})$ corresponds to the 
evaluation of $\FK{U({\bf s})}$ in $1\in U({\bf s})$.) 
If $\delta\in Der_1(\FK{U({\bf s})})$, and $V$ is an object of $\mathcal C$, the corresponding endomorphism $\delta_V\in End(V)$ is determined by
\begin{eqnarray*}
  \phi(\delta_V v)=\delta(g_{\phi v} \res{U({\bf s})} ) \quad\mb{ for all }\quad \phi\in V^{du},\;v\in V. 
\end{eqnarray*}
The Lie algebra structure on $Der_1( \FK{U({\bf s})})$ is the uniquely determined Lie algebra structure, such that the maps  
$Der_1(\FK{U({\bf s})})\to End(V)$, $V$ an object of $\mathcal C$, are morphisms of Lie algebras. 
\end{Remark}

\subsection[Prounipotent subgroups of $M$]{Prounipotent subgroups of $M$ associated to subalgebras of $\g$ which act locally nilpotent}
For a Lie algebra $\bf l$ denote the {\it lower central series} by ${\bf l}^k$, $k\in\Nn$, i.e.,
\begin{eqnarray*}
  {\bf l}^0:={\bf l}\quad \mb{ and }\quad {\bf l}^{k+1}:=[{\bf l},{\bf l}^k],\quad k\in\Nn. 
\end{eqnarray*}
We call a Lie algebra $\bf l$ {\it nilpotent} if there exists some $k\in\Nn$ such that ${\bf l}^k=\{0\}$.\vspace*{1ex} 

Now let $\n$ be a Lie subalgebra of $\g$ which acts {\it locally nilpotent} on every module $V$ contained in $\mathcal C$, i.e.,
\begin{itemize}
\item $\n$ acts locally finite on $V$, 
\item for all $v\in V$ there exists an element $k\in \Nn$ such that for all $x_0,\,x_1,\,\ldots,\,x_k\in \n$ we have
\begin{eqnarray*}
   x_0 x_1\cdots x_k v=0.
\end{eqnarray*}
\end{itemize}
Let $V$ be a module contained in $\mathcal C$. Set $V_{-1}:=\{0\}$. For $k\in \Nn$ set 
\begin{eqnarray*} 
   V_k:=\Mklz{v\in V}{ x_0 x_1\cdots x_k v=0 \;\mb{ for all }\; x_0,\,x_1,\,\ldots,\,x_k\in \n}.
\end{eqnarray*}  
We have $\n V_k\subseteq V_{k-1}$ for all $k\in\Nn$. Furthermore
\begin{eqnarray*}
   V_{k-1}\subseteq V_k \;\mb{ for all }\; k\in\Nn, \quad \mb{ and }\quad \bigcup_{i\in\Nn} V_i =V. 
\end{eqnarray*}
For $k\in\Nn$ define an ideal $I_k$ of $\n$ by
\begin{eqnarray*}
   I_k:=\Mklz{x\in \n}{xV_k=\{0\} \mb{ for all modules } V \mb{ contained in }{\mathcal C}}.
\end{eqnarray*}
The Lie algebra $\n$ acts faithfully on the modules contained in $\mathcal C$, since $\g$ does. Therefore 
\begin{eqnarray*}
   I_k\supseteq I_{k+1}\;\mb{ for all }\;k\in\Nn,\quad \mb{ and }\quad \bigcap_{k\in\Nn}I_k=\{0\}.
\end{eqnarray*}

Because of $\n^k\subseteq I_k$ the Lie algebra $\n/I_K$ is nilpotent for all $k\in\Nn$.
{\bf Assumption:} For simplicity of presentation we assume that the nilpotent
Lie algebra $\n/I_k$ is finite dimensional for all $k\in\Nn$.\vspace*{1ex}

For $k,m\in\Nn$, $k\leq m$, there are the canonical morphisms $\n / I_m \to\n / I_k$ with the obvious concatenation properties.  
The corresponding projective limit Lie algebra 
\begin{eqnarray*}
  (\n_f, (\psi_k:\n_f\to\n/I_k)_{k\in\Nn})
\end{eqnarray*}
is a pro-nilpotent Lie algebra. The Lie algebra $\n$ embeds into $\n_f$ in the obvious way.
For $k,m\in\Nn$, $k\leq m$, the morphisms $\n / I_m \to\n / I_k$ induce dual maps $(\n/I_k)^* \to (\n/I_m)^*$ with the obvious concatenation
properties. Let  
\begin{eqnarray*}
   (\n_f^{du}, (\psi_k^*: (\n/I_k)^* \to \n_f^{du})_{k\in\Nn})
\end{eqnarray*}
the corresponding direct limit linear space. Interpret $\n_f^{du}$ as a subspace
of $\n_f^*$.
For every $k\in\Nn$ we equip $\n/I_k$ with its coordinate ring $\FK{\n/I_k}$ as a finite dimensional linear space. For $k,m\in\Nn$, $k\leq m$ 
the morphisms $\n / I_m \to\n / I_k$ induce comorphisms $\FK{\n/I_k} \to\FK{\n/I_m}$ with the obvious concatenation properties. 
Let 
\begin{eqnarray*}
   (\FK{\n_f}, (\psi_k^*:\FK{\n/I_k}\to \FK{\n_f})_{k\in\Nn})
\end{eqnarray*}
be the corresponding direct limit algebra. Interpret $\FK{\n_f}$ as a coordinate ring
on $\n_f$. It is easy to check that $\FK{\n_f}$ is a symmetric algebra in $\n_f^{du}$.

For $k\in\Nn$ let $U_k$ be the unipotent linear algebraic group associated to the nilpotent Lie algebra $\n/I_k$ by the 
Campbell-Baker-Hausdorff formula, compare for example \cite{Bou}, \S9, n$\mb{}^o\,$5. Then $\n/I_k$ identifies with the Lie algebra
of $U_k$. The exponential map $\exp:\n/I_k\to U_k$ is an isomorphism of varieties. 
For $k,m\in\Nn$, $k<m$, the canonical morphism $\n / I_m \to \n / I_k$ induces a morphism $U_m\to U_k$. 
The corresponding projective limit group  
\begin{eqnarray*}
   (U_f, (\Psi_k: U_f\to U_k)_{k\in\Nn}) 
\end{eqnarray*}
is a pro-unipotent group. 
For $k,m\in\Nn$, $k\leq m$, the morphisms $U_m \to U_k$ induce comorphisms $\FK{U_k} \to\FK{U_m}$ with the obvious concatenation properties. Let
\begin{eqnarray*}
  (\FK{U_f}, (\Psi_k^*:\FK{U_k}\to \FK{U_f})_{k\in\Nn})
\end{eqnarray*}
be the corresponding direct limit algebra. Interpret $\FK{U_f}$ as a coordinate ring
on $U_f$. It is well known that $\FK{U_f}$ is in the natural way a Hopf algebra.
There is a bijective map $\exp: \n_f\to U_f$, uniquely determined by $\Psi_k(\exp(x))=\exp(\psi_k(x))$, $x\in\n_f$, $k\in\Nn$, which
we also call exponential map. It is an isomorphism of sets with coordinate rings.

Now let $V$ be a $\g$-module contained in $\mathcal C$. The Lie algebra $\n/I_k$ acts on $V_k$ in the obvious way, $k\in\Nn$. 
We get a well defined action of the Lie algebra $\n_f$ on $V$ by
\begin{eqnarray*}
   x_V v:=\psi_k(x)v \quad \mb{ where }\quad x\in \n_f\quad\mb{ and }\quad v\in V_k,\;k\in\Nn.
\end{eqnarray*}
Since $\n_f V_k\subseteq V_{k-1}$, $k\in\N$, the Lie algebra $\n_f$ acts also locally nilpotent on $V$. 
Similarly, the group $U_k=\exp(\n/I_k)$ acts on $V_k$ in the obvious way, $k\in\Nn$. We get a well defined action of the group $U_f$ on $V$ by
\begin{eqnarray*}
   u_V v:=\Psi_k(u)v \quad \mb{ where }\quad u\in U_f\quad\mb{ and }\quad v\in V_k,\;k\in\Nn.
\end{eqnarray*}
These actions are compatible with the exponential function $\exp:\n_f\to U_f$
and the exponential function on locally finite endomorphisms, i.e., 
\begin{eqnarray*}
   \exp(x)_V =\exp(x_V)\quad \mb{ for all }\quad x\in\n_f.
\end{eqnarray*}

\begin{Theorem}\label{unipotent} Let $\n$ be a Lie subalgebra of $\g$, such that $\n$ acts locally nilpotent on every module $V$ contained in 
$\mathcal C$. Let $\n/I_k$ be finite dimensional for all $k\in\Nn$. Set
\begin{eqnarray*}
     D    &:=& \Mklz{u\in U_f}{u_V \in End_{V^{du}}(V) \mb{ for all objects } V \mb{ of }{\mathcal C}},\\
  {\bf d} &:=& \Mklz{x\in \n_f}{x_V \in End_{V^{du}}(V) \mb{ for all objects } V \mb{ of }{\mathcal C}}.
\end{eqnarray*} 
Then $D$ is a submonoid of $U_f$ and $\bf d$ is a Lie subalgebra of $\n_f$. The map
\begin{eqnarray*}
  i:\; D &\to& M\\
   u  &\to& (u_V)_{V\;object\;of\;{\mathcal C}}
\end{eqnarray*}
is a morphism of monoids which is also a closed embedding of sets with coordinate
rings. If $D$ is dense in $U_f$ then the map
\begin{eqnarray*}
 Lie(i):\; {\bf d} &\to& Lie(M)\\
    x  &\to& (x_V)_{V\;object\;of\;{\mathcal C}}
\end{eqnarray*}
is an embedding of Lie algebras with image $Lie(i(U_f))$. 
\end{Theorem}

\begin{Remarks}\label{Runipotent} If $V^{du}=V^*$ for all $\g$-modules $V$ contained in $\mathcal C$, or if for all
modules $V$ contained in $\mathcal C$ the elements $x_V^{du}$, $x\in\n_f$, act
locally nilpotent on $V^{du}$, then $D=U_f$ and ${\bf d}=\n_f$. 
%
%
\end{Remarks}

\Proof The map $U_f\to\Spm\FK{U_f}$, mapping an element of $U_f$ to its evaluation homomorphism, is bijective. 
The map $\n_f\to Der_1(\FK{U_f})$, mapping an element $x\in\n_f$ to the derivation $\ti{\delta}_x$ defined by
\begin{eqnarray*}
   \ti{\delta}_x(f):= \frac{d}{dt}\res{t=0}\left(\exp^*(f)\right)(tx)= \frac{d}{dt}\res{t=0}f(\exp(tx)) 
\end{eqnarray*}
is bijective. 
Both statements are well known. They follow easily, because the map $\exp:\n_f\to U_f$ is an isomorphism of sets with coordinate rings, the algebra $\FK{\n_f}$ is a symmetric
algebra in $\n_f^{du}$, and $(\n_f^{du})^*$ identifies with the subset of $(\n_f^*)^*$ given by the evaluations in the elements of $\n_f$. (Here it is
used that $\n/I_k$ is finite dimensional for all $k\in\Nn$.)\vspace*{1ex} 

We give in two steps an isomorphism from $\FK{U_f}$ to $\FK{U(\n)}$.
The theorem then follows easily by Theorem \ref{embedding} and Remark \ref{construct}.\vspace*{1ex}

(a) With the notation of above, the coordinate ring on $U_f$ is obtained by
\begin{eqnarray*}
    \FK{U_f} =\bigcup_{k\in\Nn}\Mklz{\ti{f}\circ\Psi_k: U_f\to\F}{\ti{f}\in \FK{U_k}}.
\end{eqnarray*}
It is easy to check that $U_f$ acts on the $\g$-modules contained in $\mathcal C$ compatible with $\g$-module homomorphisms, 
compatible with tensor products, and as identity on the trivial one-dimensional $\g$-modules.

For $V$ an object of $\mathcal C$, $v\in V$ and $\phi\in V^{du}$ define the matrix coefficient $h_{\phi v}:U_f\to \F$ by 
$h_{\phi v}(u):=\phi(uv)$, where $u\in U_f$. Then
\begin{eqnarray*}
     A:=\Mklz{h_{\phi v}}{v\in V,\;\phi\in V^{du},\;V\mb{ an object of }{\mathcal C}}
\end{eqnarray*}
is a coordinate ring on $U_f$. We show that it coincides with $\FK{U_f}$.
To show $\FK{U_f}\supseteq A$, let $v\in V$, $\phi\in V^{du}$, $V$ an object of $\mathcal C$. Then there
exists an integer $k\in\Nn$ such that $v\in V_k$. Since $\n$ acts locally finite, the subspace $U(\n/I_k)v$ of
$V_k$ is finite dimensional. It is $U_k$-invariant. Therefore the matrix coefficient
$\ti{h}_{\phi v}:U_k\to \F$ defined by $\ti{h}_{\phi v}(u_k):=\phi(u_k v)$,
where $u_k\in U_k$, is an element of $\FK{U_k}$. It follows $h_{\phi v}=\ti{h}_{\phi v}\circ\Psi_k \in \FK{U_f}$.

To show $\FK{U_f}\subseteq A$ fix $k\in\Nn$. The Lie algebra $\n/I_k$ acts faithfully on the modules $V_k$, 
$V$ an object of $\mathcal C$. Therefore it also acts faithfully on the finite-dimensional modules
\begin{eqnarray*}
   U(\n/I_k)v \quad \mb{ where }\quad v\in V_k,\; V\mb{ an object of }{\mathcal C}.
\end{eqnarray*}
Since $\n/I_k$ is finite dimensional it follows easily that we may choose
finitely many modules $U(\n/I_k)v_1$, \ldots, $U(\n/I_k)v_m$, such that $\n/I_k$ acts
faithfully on $S:=U(\n/I_k)v_1\oplus \cdots \oplus U(\n/I_k)v_m$. Also $U_k$
acts faithfully on $S$.
We have $S\subseteq W_k$ for some $\g$-module $W$ of $\mathcal C$. Because $S$
is finite-dimensional, the restriction of $W^{du}$ onto $S$ coincides with $S^*$
With the same notation for the matrix coefficients of $U_k$ as above, the coordinate ring $\FK{U_k}$ is generated by the functions 
$\ti{h}_{\phi v}$, $\phi\in S^*$, and $v\in S$. (Use for example Theorem 1.3
in combination with Theorem 1.1 of Chapter VIII of the book \cite{Ho}.)Now the inclusion follows,
because $A$ is an algebra, and $h_{\phi v}=\ti{h}_{\phi  v}\circ\Psi_k \in\FK{U_f}$ for all $\phi\in S^*$, and $v\in S$. \vspace*{1ex}

(b) For $u\in U_f$ denote by $l_u$ the left multiplication of $U_f$ by $u$,
and by $l_u^*$ its comorphism. It is easy to check that for every $x\in\n$ we
get a well defined derivation of $\FK{U_f}$ by $(x_\ro h)(u):=\ti{\delta}_x(l_u^* h)$. In this way we get
an action of $\n$ on $\FK{U_f}$, which extends to an action of $U(\n)$ on
$\FK{U_f}$. It satisfies $(x_\ro h_{\phi v})=h_{\phi \,xv}$, where $x\in U(\n)$,
$v\in V$, $\phi\in V^{du}$, and $V$ an object of $\mathcal C$. 
It follows in the same way as in the proof of Theorem \ref{UnivEin} (b) that we get a surjective homomorphism of algebras 
$\Psi:\FK{U_f}\to\FK{U(\n)}$ by $\Psi(h)(x):=(x_\ro h)(1)$ which satisfies
$\Psi(h_{\phi v})=g_{\phi v}\res{U({\bf n})}$. 
To show that it is injective let $V$ be an object of $\mathcal C$, let $v\in V$ and $\phi\in V^{du}$ such that
$\phi(U(\n)v)=0$. Since $U(\n)v$ is $\n_f$-invariant and therefore also $U_f$-invariant we get $\phi(U_f v)=0$.
\qed

\subsection[Generalized toric submonoids of $M$]{Generalized toric submonoids of $M$ associated to subalgebras of $\g$ which act diagonalizable}
Let $(A,+)$ be an abelian monoid. The set of homomorphisms of monoids
\begin{eqnarray*}
  \widetilde{A}:= \mb{Hom}\left((A,+), (\F,\cdot)\right)
\end{eqnarray*}
gets the structure of an abelian monoid by multiplying the homomorphisms pointwise. The monoid algebra $\FK{A}$ can be identified 
with a coordinate ring $\FK{\widetilde{A}}$ on $\widetilde{A}$, identifying $\sum_{a\in A} c_a a\in\FK{A}$ (which means $c_a\in\F$, $c_a\neq 0$ only for finitely 
many $a\in A$) with the function
\begin{eqnarray*}
   (\sum_{a\in A} c_a a)(\al):=\sum_{a\in A} c_a \al(a),\quad \al\in \widetilde{A}.
\end{eqnarray*}
Equip the $\F$-linear space of homomorphisms of monoids 
\begin{eqnarray*}
   Lie(\widetilde{A}):= \mb{Hom}\left((A,+), (\F,+)\right).
\end{eqnarray*}
with the structure of an abelian Lie algebra.
\begin{Theorem}\label{Tor1} Let $\h$ be an abelian subalgebra of $\g$, such that every $\g$-module $V$ contained in $\mathcal C$ has a weight 
space decomposition $V=\bigoplus_{\la\in {\bf h}^*} V_\la$. Denote by $P(V)\subseteq \h^*$ the set of weights of $V$. Then
\begin{eqnarray*}
  A:=\bigcup_{V\; an\; object\;of\;{\mathcal C}} P(V)
\end{eqnarray*}
is a submonoid of $(\h^*,+)$. Its $\F$-linear span in $\h^*$ is point separating on $\h$. Furthermore:  
\begin{itemize}
\item[(a)] Let $D$ be the set of elements $\al\in\widetilde{A}$ such that for every $\g$-module $V$ contained in $\mathcal C$ the endomorphism  
$\al_V\in  End(V)$  defined by
\begin{eqnarray*}
     \al_V v_\la := \al(\la)v_\la,\quad v_\la\in V_\la,\quad\la\in P(V),
\end{eqnarray*}
is already in $End_{V^{du}}(V)$. Then $D$ is a submonoid of $\widetilde{A}$ and the map  
\begin{eqnarray*}
 t:\; D &\to& M\\
   \al &\mapsto & (\al_V)_{V\;object\;of\;{\mathcal C}}
\end{eqnarray*}
is an morphism of monoids which is also a closed embedding of sets with coordinate rings.
\item[(b)] Let $\bf d$ be the set of elements $\beta\in Lie(\widetilde{A})$ such that for every $\g$-module $V$ contained in $\mathcal C$ the 
endomorphism  $\beta_V\in End(V)$  defined by
\begin{eqnarray*}
     \beta_V v_\la := \beta(\la)v_\la,\quad v_\la\in V_\la,\quad\la\in P(V),
\end{eqnarray*}
is already in $End_{V^{du}}(V)$. Then $\bf d$ is a Lie subalgebra of
$Lie(\widetilde{A})$. If $D$ is dense in $\widetilde{A}$ then the map 
\begin{eqnarray*}
 Lie(t):\;Lie({\bf d}) &\to& Lie(M)\\
   \beta &\mapsto & (\beta_V)_{V\;object\;of\;{\mathcal C}}
\end{eqnarray*}
is an embedding of Lie algebras with image $Lie(t(D))$. 
\end{itemize}
\end{Theorem}
\begin{Remark}\label{RTor} If $V^{du}=V^*$ for all $\g$-modules $V$ contained in $\mathcal C$, or if 
  or if $V^{du}$ has a weight space decomposition with respect to $\h$ for all
  $\g$-modules $V$ contained in $\mathcal C$, then $D=\widetilde{A}$ and ${\bf d}=Lie(\widetilde{A})$.
\end{Remark}

\Proof Since $P(V_0)=\{0\}$ and $P(V\otimes W)=P(V)+P(W)$ for all objects $V$, $W$ of $\mathcal C$, the set $\widetilde{A}$ is a submonoid 
of $(\h^*,+)$. Let $h\in\h$ such that $\phi(h)=0$ for all $\phi\in A$. Then $h$ acts trivially on every $\g$-module contained in $\mathcal C$. 
Therefore $h=0$.

Let $V$ be an object of $\mathcal C$, $v\in V\setminus\{0\}$, and $\phi\in V^{du}$. For
$\la\in P(V)$ denote by $v_\la$ the $V_\la$-weight space component of $V$, and
by by $supp(v)$ the set of weights for which $v_\la\neq 0$. It is easy to
check that the matrix coefficient of $\phi$ and $v$ on $U(\h)$ is given by
\begin{eqnarray}\label{mctoricexplicit}
  g_{\phi v}\res{U({\bf h})}=\sum_{\la\in supp(v)} \phi(v_\la)\la,
\end{eqnarray}
where $\la\in U(\h)^*$ is determined by 
\begin{eqnarray*}
 \la(1)=1 \quad\mb{ and }\quad \la(h_1\cdots h_n)=\la(h_1)\cdots \la(h_n) \quad\mb{ where }\quad h_1,\ldots,h_n\in\h,\;n\in \N.
\end{eqnarray*}
It is easy to check that $\FK{U(\h)}=\bigoplus_{\la\in
  A} \F \la$ identifies with the monoid algebra $\FK{A}$. Restricting the elements of
$\Spm\FK{A}$ to $A\subseteq \FK{A}$, the set $\Spm\FK{A}$ identifies with
$\widetilde{A}$. Similarly, restricting the elements of
$Der_1\FK{A}$ to $A\subseteq \FK{A}$, the set $Der_1\FK{A}$ identifies with
$Lie(\widetilde{A})$. The theorem now follows by Theorem \ref{embedding}
and Remark \ref{construct}. (Use equation (\ref{mctoricexplicit})). 
\qed

In the following Theorem we state an description of the monoid $\ti{A}$ in the important case where $A$ is contained in a lattice of finite rank. 
We call such an monoid $\ti{A}$ a {\it generalized toric monoid} since it has similar properties as
an algebraic toric monoid.

We call a submonoid $A$ of a free $\Z$-module $L$ saturated, if for all $x\in L$ and $n\in \N$ from $nx\in A$ follows $x\in A$. 
A face $F$ of $A$ is a submonoid of $A$ such that $A\setminus F + A\subseteq A\setminus F$. We denote by $Fa(A)$ the set of faces of $A$. As
described in Theorem 1.2 and 1.3 of \cite{M1} the faces of $A$ correspond bijectively to the faces of the convex cone generated by $A$ in the real
vector space $L\otimes_\Z\R$.

In the subsection '1.4 Generalized toric varieties' of \cite{M1}, we gave
the following Theorem which fits very well into the situation which we
consider here.

\begin{Theorem}\label{Tor2} Let $A$ be a saturated submonoid of a free
  $\Z$-module of finite rank. Let $Fa(A)$ be its set of faces. For $F\in\Fa{A}$ set
\begin{eqnarray*} 
     T(F)  &:=&  \Mklz{\al\in\widetilde{A}\,}{\,\al^{-1}(\F^\times)=F }\;\;.
\end{eqnarray*} 
Let $e(F)\in T(F)$ be the element given by
\begin{eqnarray*}
 e(F) \la \;:=\; \left\{\begin{array}{cl}
                  1 & \mb{ if }\;\la\in F \\ 
                  0 & \mb{ if }\;\la\in A\setminus F
                  \end{array}  \right.\quad .
\end{eqnarray*}
Then we have:
\begin{itemize}

\item[(1)] $T(A)$ is the unit group of $\widetilde{A}$, and $E:=\Mklz{e(F)}{F\in\Fa{A}}$ is the set of idempotents of
$\widetilde{A}$. Furthermore, $\widetilde{A}= T(A)\, E = E\, T(A) $.
\item[(2)] $T(F)$ is a subgroup of $\widetilde{A}$ with unit $e(F)$, isomorphic to the torus $\widetilde{F-F}$, an isomorphism
$\Phi_F:\,\widetilde{F-F}\to T(F)$ given by
\begin{eqnarray*}
    \Phi_F(\al)(\la)\;:=\;\left\{ \begin{array}{cl}
                                  \al(\la) & \la\in F \\
                                     0   & \la\in A\setminus F
                                 \end{array}\right. &\;\;,\;\; & \al\in \widetilde{F-F}\;\;.       
\end{eqnarray*}
Let the torus $T(A)$ act on $\widetilde{A}$ by left multiplication. The tori $T(F)$, $F\in\Fa{A}$, are the $T(A)$-orbits of $\widetilde{A}$.
\item[(3)] For $\la\in A$ denote by $D(\la)$ the principal open set $\{\al\in\widetilde{A}\,|\,\al(\la)\neq 0\}$. We have
\begin{eqnarray*}
  \overline{T(F)} &=&   \bigcup_{G\in\Fa{F}} T(G)\;\;, \\
         D(\la)    &=&   \bigcup_{G\in \Fa{A}\,,\,G\supseteq F } T(G)  \quad \mb{ where }\quad  \la\in ri\,F\;\;.
\end{eqnarray*}
In particular, $T(A)$ is principal open and dense in $\widetilde{A}$.
\end{itemize}
\end{Theorem}
\subsection{The case where $\g$ is generated by integrable locally finite elements}

From the embedding theorem \ref{embedding} follows easily: 

\begin{Corollary}\label{integrablelocallyfinite} Let $e\in \g$ act locally finite on every $\g$-module contained in $\mathcal C$. 
Suppose there exists a dense submonoid $D$ of $M(\F e)$ such that for all $d\in D$ 
and all $\g$-modules $V$ contained in  $\mathcal C$  we have $d_V\in End_{V^{du}}(V)$. 
\begin{itemize}
 \item[(a)] Then the map 
\begin{eqnarray*}
  i:  D  &\to     &  M\\
   d &\mapsto & (d_V)_{V\;object \;of\;{\mathcal C}} 
\end{eqnarray*}
is a morphism of monoids which is also an embedding of sets with coordinate rings.
\item[(b)] We have $\F e \subseteq Lie(i(D))\subseteq Lie(M)$.
\end{itemize}  
\end{Corollary}

\begin{Definition} We call an element $e\in\g$ as in the last corollary an integrable locally finite element of $\g$ (with respect to $\mathcal C$,
  $\mathcal C^{du}$). We call the map $i:D\to M$ a parametrization
  corresponding to $e$. If $D=\Mklz{d\in M(\F e)}{d_V\in End_{V^{du}}(V)}$ then
  we call $i:D\to M$ the maximal parametrization corresponding to $e$.
\end{Definition}

There are two sorts of integrable locally finite elements which are particularly important. These are described in the next two theorems. 
It would be possible to use Theorem \ref{unipotent}, Theorem \ref{Tor1}, and Theorem
\ref{Tor2} for their proofs. We give here direct proofs. Some of the results
obtained in these proofs are needed in \cite{M4}.

\begin{Theorem} \label{one1} Let $x\in \g\setminus\{0\}$ such that for every $\g$-module $V$ contained in $\mathcal C$ 
\begin{itemize}
\item[(1)] $x_V$ is locally nilpotent,
\item[(2)] $\exp(t x_V)\in End_{V^{du}}(V)$ for all $t\in\F$.
\end{itemize}
Then $x$ is an integrable locally finite element of $\g$. Furthermore:
\begin{itemize}
\item[(a)] For every $t\in \F$ there exists an element $\exp(tx)\in M$ which acts on every $\g$-module $V$ contained 
in $\mathcal C$ by $\exp(t x_V)$.
\item[(b)] Equip $\F$ with the coordinate ring of polynomial functions $\FK{t}$. Then the map
\begin{eqnarray*}
  u_x:\; (\F,+)   &\to      & \;\;\;M\\
          t \;\;\;\;  & \mapsto & \exp(tx) 
\end{eqnarray*}
is a morphism of monoids which is also a closed embedding of sets with
coordinate rings. We denote its image by $U_x$.
\item[(c)] We have $Lie(U_x)=\F x$.
\end{itemize}
\end{Theorem}
\begin{Remarks}\label{Rone1} (1) The condition (2) of the proposition is satisfied if $V^{du}=V^*$ or if 
$(x_V)^{du}$ is a locally nilpotent endomorphism of $V^{du}$.

(2) Let $x\in\g$ satisfy the conditions (1) and (2) of the proposition. For any $m\in M^\times$ also $mx\in \g$ satisfies the conditions of 
the proposition and
\begin{eqnarray*}
    m\exp(x) m^{-1} =\exp(mx).
\end{eqnarray*}  
\end{Remarks}

\Proof By the Poincare-Birkhoff-Witt theorem $U(\F x)=\oplus_{n\in\Nn}\F x^n$. Define
$\tau\in U(\F x)^*$ by $\tau(\frac{x^m}{m!})=\delta_{1 m}$, $m\in\Nn$. It is
easy to check that
\begin{eqnarray}\label{taudu}
    \tau^n(\frac{x^m}{m!})=\delta_{n m}\quad \mb{ where }\quad  n,\,m\in\Nn. 
\end{eqnarray}
We show
\begin{eqnarray}\label{mcU(Fx)}
   \FK{U(\F x)}=\bigoplus_{n\in\Nn}\F \tau^n .
\end{eqnarray}
To show the inclusion ``$\subseteq $'' let $V$ be a $\g$-module contained in $\mathcal C$, let $v\in V$ and $\phi\in V^{du}$. Since $x$ acts
locally nilpotent on $V$ there exists an element $n\in\Nn$ such that
$x^{n+1}v=0$. With (\ref{taudu}) follows
\begin{eqnarray}\label{mcexplicit}
      g_{\phi v}\res{U(\F x)}=  \sum_{k=0}^n \phi(\frac{x^k}{k!} v) \tau^k.
\end{eqnarray} 
To show the inclusion ``$\supseteq $'' we have to show $\tau\in \FK{U(\F x)}$. Since $\g$ acts faithfully 
on the $\g$-modules contained in $\mathcal C $, there exists a $\g$-module $V$ contained in $\mathcal C$ such that $x_V\neq 0$. 
Since $x_V$ is a locally nilpotent endomorphism there exists an element $v\in V$
such that $v\neq xv\neq 0$ and $x^2 v=0$.
Since $V^{du}$ separates the points of $V$ there exists an element $\psi\in V^{du}$ 
such that $\psi(v)=0$ and $\psi(xv)=1$. By equation (\ref{mcexplicit}) we get
$g_{\phi v}\res{U(\F x)}=\tau$.

By (\ref{mcU(Fx)}) the algebra $\FK{U(\F x)}$ coincides with the polynomial
algebra $\FK{\tau}$. The map $\F\to\Spm\FK{\tau}$, which maps $t\in\F$ to the
evaluation morphism in $t$, is bijective. The map $\F\to Der_1\FK{\tau}$, which maps $t\in\F$ to $t\frac{d}{d\tau}\res{\tau=0}$, is
bijective. The Theorem now follows by Theorem \ref{embedding} and Remark
\ref{construct}. (Use equation (\ref{mcexplicit})).
\qed

Let $W$ be a $\F$-linear space. We say $\al\in End(W)$ is diagonalizable with
integer eigenvalues if there exists a base $(w_j)_{j\in J}$ of $W$ such that
\begin{eqnarray*}
   \al w_j=z_j w_j \quad\mb{ with }\quad z_j\in\Z\quad\mb{ for all }\quad j\in J.
\end{eqnarray*}
For $s\in\F^\times$ we define the endomorphism 
\begin{eqnarray*}
  s^\al \in End(W) \quad \mb{ by }\quad   s^\al w_j := s^{z_j}w_j \quad \mb{ for all }\quad j\in J.
\end{eqnarray*}

\begin{Theorem}\label{one2}
Let $h\in\g$ such that for every module $V$ contained in $\mathcal C$ 
\begin{itemize}
\item[(1)] $h_V$ is diagonalizable with integer eigenvalues,
\item[(2)] $s^{h_V}\in End_{V^{du}}(V)$ for every $s\in\F^\times$. 
\end{itemize}
Then $h$ is an integrable locally finite element of $\g$. Furthermore:
\begin{itemize}
\item[(a)] For every $s\in\F^\times$ there exists an element $s^h\in M$ which acts on every module $V$ contained in 
$\mathcal C$ by $s^{h_V}$.
\item[(b)] The set $EV$ of all eigenvalues of $h_V$ for all objects $V$ of $\mathcal C$ is a submonoid of $\Z$. 
Equip $\F^\times$ with the coordinate ring of Laurent polynomial functions $\FK{s, s^{-1}}$. Then
\begin{eqnarray*}
  t_h:\;(\F^\times ,\cdot)  &\to      & M\\
   s\;\;\;\;  & \mapsto &  s^h
\end{eqnarray*}
is a morphism of monoids which is also a morphism of sets with coordinate rings. We denote its image by $T_h$.
The image of the comorphism $t_h^*$ is given by
\begin{eqnarray*}
   t_h^*(\FK{M})=\bigoplus_{b\in EV} \F \,s^b
\end{eqnarray*}
\item[(c)] We have $Lie(T_h)=\F h$.  
\end{itemize}
\end{Theorem}
\begin{Remarks}\label{Rone2} (1) The condition (2) of the proposition is satisfied if $V^{du}=V^*$ or if 
$(h_V)^{du}$ is a diagonalizable endomorphism of $V^{du}$.

(2) Suppose that $h$ acts by positive and negative eigenvalues on the objects of $\mathcal C$. Let $a$ be the smallest absolute value of the non-zero 
eigenvalues of $h_V$ for all objects $V$ of $\mathcal C$. Then the kernel of the morphism $t_h$ consists of the group of $a$-th roots of unity. 
The image of the comorphism $t_h^*$ is given by $t_h^*(\FK{M})=\FK{t^a,t^{-a}}$.  

(3) Let $h\in\g$ satisfy the conditions (1) and (2) of the proposition. For any $m\in M^\times$ also $mh\in \g$ satisfies the conditions of 
the proposition and
\begin{eqnarray*}
    m s^h m^{-1} =s^{mh} \quad\mb{ for all } \quad s\in \F^\times.
\end{eqnarray*}  
\end{Remarks}

\Proof By the Poincare-Birkhoff-Witt theorem $U(\F h)=\oplus_{n\in\Nn}\F
h^n$. For $n\in\Z$ define $\exp(n\tau)\in U(\F h)^*$ by $\exp(n\tau)(\frac{h^m}{m!})=n^m$, $m\in\Nn$. It is
not difficult to check that the elements $\exp(n\tau)$, $n\in\Z$ are linearly
independent, and $\exp(n\tau)=\exp(\tau)^{-n}$ for all $n\in\Z$.
Let $V$ be a $\g$-module contained in $\mathcal C$, let $v\in V\setminus\{0\}$ and $\phi\in V^{du}$. 
Write $v$ as a finite sum $v=\sum_{i}v_i$, where $v_i$ is an eigenvector of $h$ to the eigenvalue $n_i\in EV$. It is easy to check that
\begin{eqnarray*}
   g_{\phi v}\res{U(\F h)}=\sum_i \phi(v_i) \exp(n_i \tau).
\end{eqnarray*}
It follows
\begin{eqnarray}\label{mcU(Fh)}
   \FK{U(\F h)}=\bigoplus_{n\in EV}\F \exp(\tau)^n .
\end{eqnarray}

An element $\al\in M(\F h)=\Spm\FK{U(\F h )}$ acts on an eigenvector $v\in V$ to the eigenvalue $n\in EV$, $V$ an object of $\mathcal C$, by
\begin{eqnarray*}
    \al_V v=\al (\exp(\tau)^n)v.
\end{eqnarray*} 
For $s\in \F^\times$ denote by $\al_s\in M(\F h)$ the 
homomorphism $\al_s(\exp(\tau)^n):=s^n$, $n\in EV$. By our assumptions
\begin{eqnarray*}
  D:=\Mklz{\al_s}{s\in\F^\times}\subseteq \Mklz{\al\in M(\F h)}{\al_V\in  End_{V^{du}}(V) \mb{ for all objects } V \mb{ of }{\mathcal C}}.
\end{eqnarray*}
It is easy to check that $D$ is dense in $M(\F h)$. It follows by definition that $h$ is a integrable locally finite element of $\g$. 
From part (a) of Corollary \ref{integrablelocallyfinite} follows that we get a morphisms of monoids $i:D\to M$ by
$i(\al_s)=((\al_s)_V)_{V\;an\;object\;of\; {\mathcal C}}$, which is also an embedding of sets with coordinate rings. In particular
this shows part (a) of the theorem.

(b) Equip $\F^\times$ with its coordinate ring of Laurent polynomial functions $\FK{s,s^{-1}}$. Define the map $p:\F^\times \to D$ by
$p(s):=\al_s$. Then $p$ is a surjective morphism of monoids which is also a
morphism of sets with coordinate rings. It is easy to check 
that $p^*(\FK{D})=\bigoplus_{b\in EV}\F s^b$.
Now part (b) of the theorem follows since the map $t_h:F^\times\to M$ is the concatenation of $p:\F^\times \to D$ and $i:D\to M$.

(c) From part (b) of Corollary \ref{integrablelocallyfinite} follows $\F h\subseteq Lie(T_h)$. The derivation $\ti{\delta}_h$ of
$\FK{D}$ in $1\in D$ which corresponds to $h$, is given by $\ti{\delta}_h(\exp(\tau)^b)=b$, $b\in EV$. 
To show $\F h=Lie(T_h)$ it is sufficient to show that every derivation $\ti{\delta}$ of
$\FK{D}$ in $1\in D$ is a scalar multiple of $\ti{\delta}_h$.
Let $a\in EV\setminus\{0\}$. Since $(\exp(\tau)^b)^a=\exp(ab\tau)=(\exp(\tau)^a)^b$ for all $b\in EV$ we find
\begin{eqnarray*}
   a \ti{\delta}(\exp(\tau)^b)=b \ti{\delta}(\exp(\tau)^a).
\end{eqnarray*}
It follows $\ti{\delta} = \frac{\ti{\delta}(\exp(\tau)^a)}{a} \ti{\delta}_h $.
\qed

There are very important examples of Lie algebras $\g$ and categories $\mathcal C$, $\mathcal C^{du}$, where $\g$ is generated as Lie algebra 
by integrable locally finite elements. Certain classes of examples are investigated in \cite{M4}. With the methods developed 
so far we now can show the following important theorem.
\begin{Theorem}\label{glocalfinite} Let the Lie algebra $\g$ be generated by a set
  $E$ of integrable locally finite elements. Then $\mathcal C$, ${\mathcal C}^{du}$ is good for integrating $\g$. 

For $e\in E$ let $i_e: D_e\to M$ be a corresponding parametrization. Set $V_e:=i(D_e)$. Then 
\begin{eqnarray*}
  M_E := \bigcup_{m\in\N}\quad \bigcup_{e_1,\,e_2,\,\ldots,\,e_m\in E\setminus\{0\}} V_{e_1}V_{e_2}\cdots V_{e_m}
\end{eqnarray*} 
is a dense submonoid of $M$. 

If for every $e\in E$ there exists a parametrization $i_e: D_e\to M$ such that
$D_e$ is a group then the pair $\mathcal C$, ${\mathcal C}^{du}$ is very good for integrating $\g$.
\end{Theorem}
\begin{Remarks} (1) The monoid $M_E$ can be constructed in a similar way as the Kac-Moody group in \cite{KP1}, or as 
the group associated to an integrable Lie algebra in \cite{K1}, \S 1.5, and \S
1.8, as a quotient, defined in an appropriate way, of the free product of monoids $\bigstar_{e\in E\setminus\{0\}} V_e$.

(2) If $\g$ is generated by integrable locally finite elements as used in Theorem \ref{one1} and Theorem \ref{one2}, 
then the pair $\mathcal C$, ${\mathcal C}^{du}$ is very good for integrating $\g$.

Let ${\mathcal C}^{du}={\mathcal C}^{full}$. If $\g$ is generated by
integrable locally finite elements then the pair $\mathcal C$, 
${\mathcal C}^{full}$ is very good for integrating $\g$. (This can be seen as
follows. Let $e\in E$. By Theorem \ref{fdgood} the unit group $M(\F e)^\times$ is dense in $M(\F e)$. Since
$i:\,M(\F e)\to M$ is a parametrization, also the restriction $i:\,M(\F e)^\times\to M$ is a parametrization.)
\end{Remarks}
%
%
%
%
%
%

\Proof By Corollary \ref{integrablelocallyfinite} we have $E\subseteq Lie(M)$. Since $E$ generates the Lie algebra $\g$ also 
$\g\subseteq Lie(M)$. 

Let $V$ be a $\g$-module contained in $\mathcal C$, let $v\in V$, $\phi\in V^{du}$ such that $f_{\phi v}$ vanishes on $M_E$. 
Let $e\in E\setminus\{0\}$. Then
\begin{eqnarray*}
     f_{g^{du}\phi\,v}  \quad \mb{ vanishes on } V_e \quad \mb{ for all }\quad g\in M_E.
\end{eqnarray*}
By  Corollary \ref{integrablelocallyfinite} we have $e\in Lie(V_e)$, from which follows
\begin{eqnarray*}  
  0=\delta_e \left(f_{g^{du}\phi\,v}\right)= \phi(gev) \quad \mb{ for all }\quad g\in M_E.
\end{eqnarray*}
Therefore $f_{\phi \,ev}$ vanishes on $M_E$. Repeating this argument and evaluating in $1\in M$, we find
\begin{eqnarray*}
   \phi (e_1\cdots e_m v)=0 \quad \mb{ for all } \quad e_1,\,\ldots,\, e_m\in E\setminus\{0\},\;m\in \N.
\end{eqnarray*}
Since $E$ generates the Lie algebra $\g$ this equation is also valid for all $e_1,\,\ldots,\, e_m\in \g$, $m\in\N$. Since $1\in M_E$ 
also $\phi(v)=0$. With Theorem \ref{UnivEin} (b) follows $f_{\phi v}=0$.

Now suppose that for every $e\in E$ there exists a parametrization $i_e: D_e\to M$ such that $D_e$ is a group. Then the corresponding 
dense subset $M_E$ is a subgroup of $M^\times$. By Proposition \ref{Defrel} the pair $\mathcal C$, ${\mathcal C}^{du}$ is 
very good for integrating $\g$
\qed

As already stated in Remark \ref{RemInv} above, if the Lie algebra $\g$ is
generated by integrable locally finite elements then there is an easy relation 
between the $\g$ and $M$-invariant subspaces of a $\g$-module contained in $\mathcal C$:

\begin{Theorem}\label{sub3}  Let the Lie algebra $\g$ be generated by integrable locally finite elements. Let $U$ be a linear subspace of a 
$\g$-module $V$ contained in $\mathcal C$. Then $U$ is $\g$-invariant if and only if $U$ is $M$-invariant.
\end{Theorem}

\Proof  Because of Proposition \ref{sub} we only have to show that if $U$ is
$M$-invariant it is also $\g$-invariant. It is sufficient to show that $U$ is
invariant under every integrable locally finite element. Let $e$ be an integrable locally finite
element of $\g$ with corresponding parametrization $i:D\to M$.
Since $i(D)$ acts locally finite on $V$ it also acts locally finite on
$U$. Let $W$ be a finite-dimensional $i(D)$-invariant subspace of $U$. 
Since $e\in Lie(i(D))\subseteq Lie(M)$ there exists a derivation $\ti{\delta}_e\in
Der_1\FK{i(D)}$ such that 
\begin{eqnarray*}
\ti{\delta}_e(f_{\phi w}\res{i(D)})=\delta_e(f_{\phi w})=\phi (e_V w) \quad \mb{ for all }\quad\phi\in V^{du},\;w\in W.
\end{eqnarray*}
Let $w_1,\,\ldots,\,w_n$ a base of $W$. Let $\phi_1,\,\ldots,\,\phi_n$ be
elements of $V^{du}$ which restrict on $W$ to the dual base. Define a map
$\ti{e}_V:W\to W$ by 
\begin{eqnarray*}
  \ti{e}_V w:= \sum_{i=1}^n \ti{\delta}_e(f_{\phi_i w}\res{i(D)})\, w_i\quad \mb{ for all } \quad w\in W.
\end{eqnarray*}
Because of $\sum_{i=1}^n \phi(w_i)\phi_i\res{W}= \phi\res{W}$ we find
\begin{eqnarray*}
   \phi(\ti{e}_V w)=\sum_{i=1}^n \ti{\delta}_e(f_{\phi_i w}\res{i(D)})\phi(w_i)=\ti{\delta}_e(f_{\sum_{i=1}^n \phi(w_i)\phi_i w}\res{i(D)})
   =\ti{\delta}_e(f_{\phi w}\res{i(D)})=\phi(e_V w)
\end{eqnarray*}
for all $w\in W$. It follows $e_V w=\ti{e}_V w\in W$ for all $w\in W$.
\qed
%
%
%

%
%
%
\section{The algebraic geometry of the Tannaka reconstruction\label{AGeo}}
%
%
%
%
%

In this section we interprete any Tannaka monoid $M$ associated to pair of
categories $\mathcal C$, ${\mathcal C}^{du}$ algebraic geometrically, as weak algebraic monoid with Lie
algebra $Lie(M)$, acting on the objects of $\mathcal C$. 
We show that the action of the Lie algebra $Lie(M)$ on any
object of $\mathcal C$, as introduced in the last section, is the differentiated action of the monoid $M$.

The algebraic geometry has already been developed and used in \cite{M1}. The
starting point have been 
the topologized coordinate rings of the integrable highest weight modules
$L(\La)$ of a symmetrizable Kac-Moody algebra, and of the Kostant cones
${\mathcal V}_\La\subseteq L(\La)$ given in \cite{KP2}, Section 3A.
We state here the facts, which we use in this paper and in \cite{M4}.
In difference to \cite{M1} we introduce first, as an intermediate step in the
definition of varieties and weak algebraic monoids, F-varieties and weak
F-algebraic monoids. Also we have no need to define and work with non-closed
varieties. Because of both it easier to follow the definitions and constructions.
The definition of a weak algebraic monoid is slightly more restrictive as in
\cite{M1}, Chapter 3, Section 2, to fit optimal to the Tannaka monoids. We also introduce weak algebraic groups, and
actions of weak algebraic monoids and groups. We give some remarks
concerning the relation to affine ind-varieties.

\subsection{F-varieties and weak F-algebraic monoids and groups\label{F-var}}
Recall the category of sets with coordinate rings from Section \ref{Pre}. Recall that for a set with coordinate ring $(A,\FK{A})$ 
we denote by $\Spm\FK{A}$ the set of $\F$-valued points of $\FK{A}$, i.e., the set of homomorphism of algebras from $\FK{A}$ to $\F$. 
We get an injective map
\begin{eqnarray*}
   A \to \Spm\FK{A}
\end{eqnarray*}
by assigning $a\in A$ its evaluation homomorphism $\al_a$, defined by $\al_a(f):=f(a)$ for all $f\in\FK{A}$. 
To cut short our notation we identify $A$ with its image in $\Spm\FK{A}$.\vspace*{1ex}

Let $(A,\FK{A})$ be a set with coordinate ring. A non-empty subset $\Fi{A}$ of
ideals of $\FK{A}$ is called a {\it filter of ideals} if the following holds: If $I,\,J\in\Fi{A}$ then also $I\cap J\in
\Fi{A}$. Furthermore, if $I\in\Fi{A}$ and $J$ is an ideal of $\FK{A}$ such that $I\subseteq J$, then also $J\in \Fi{A}$.
The set of ideals of $\FK{A}$ itself is a filter, which we call the {\it discrete filter}.

A non-empty subset ${\mathcal B}_A$ of the filter $\Fi{A}$ is called {\it filterbase}, if for
all $I\in\Fi{A}$ there exists a $B\in {\mathcal B}_A$ such that $B\subseteq I$.
A non-empty set ${\mathcal B}_A$ of ideals of $\FK{A}$ is a filterbase for some
filter of ideals, if and only if for all $B_1,\,B_2 \in{\mathcal B}_A$ there
exists an element $B_3\in{\mathcal B}_A$ such that $B_3\subseteq B_1\cap B_2$.\vspace*{1ex}

{\bf The category of F-varieties:} A {\it F-variety} consists of a non-empty set $A$ equipped with a coordinate ring $\FK{A}$,
 and a filter of ideals $\Fi{A}$ of $\FK{A}$, such that
\begin{eqnarray*}
       A = \Mklz{\al\in\Spm\FK{A} }{ \mb{Kernel }\al \in \Fi{A}}.
\end{eqnarray*}
A {\it morphism} of the F-variety $(A,\FK{A},\Fi{A})$ to the F-variety $(A,\FK{A},\Fi{B})$ consists of a map  $\phi:A\to B$, 
such that the comorphism $\phi^*:\FK{B}\to\FK{A}$ exists, and satisfies $(\phi^*)^{-1}\left(\Fi{A}\right)\subseteq \Fi{B}$.

Define the {\it tangent space} of a F-variety $(A,\FK{A},\Fi{A})$ at a point $a\in A$ by
\begin{eqnarray*}
      T_a A := \Mklz{\delta\in\mb{Der}_a\FK{A} }{I\subseteq 
     \mb{Kernel }\delta \;\;\mb{ for some } \;\;I\in \Fi{A}}.
\end{eqnarray*}
%
%
If $(A,\FK{A},\Fi{A})$, $(B,\FK{B},\Fi{B})$ are F-varieties, $a\in A$, and if $\phi: A\to B$ is a morphism, we get
a linear map, the {\it tangent map} at $a$, by:
\begin{eqnarray*}
 T_a\phi:\;\;T_a A &\to & T_{\phi(a)}B \\
  \delta\quad  &\mapsto & \delta\circ\phi^*
\end{eqnarray*}

{\bf Substructures:} Let $(B,\FK{B},\Fi{B})$ be a F-variety. Let $A\subseteq B$ be Zariski closed,
and $a\in A$. The {\it F-subvariety} 
$(A,\FK{A},\Fi{A})$ is defined by restricting the functions of 
$\FK{B}$ onto $A$, i.e.,
\begin{eqnarray*}
   \FK{A} := \FK{B}\res{A} \quad \mb{ and }\quad \Fi{A} := \Mklz{ I\res{A} }{ I\in \Fi{B} }.
\end{eqnarray*}
Denote by $I_B(A)$ the vanishing ideal of $A$ in $\FK{B}$. The inclusion map
$\iota: A\to B$ is a morphism. Its tangent map 
$T_a\iota:\,T_a A\to T_a B$ is injective with image
\begin{eqnarray*} 
T_a\iota(T_a A) = \Mklz{\delta\in T_a B}{\mb{ Kernel }\delta\supseteq I_B(A)} . 
\end{eqnarray*}

 {\bf Products:} Let $(C,\FK{C},\Fi{C})$, $(D,\FK{D},\Fi{D})$ be F-varieties. Equip the coordinate ring $\FK{C\times D}\cong\FK{C}\otimes\FK{D}$
 of the product $C\times D$ with the filter of ideals determined by the filterbase 
\begin{eqnarray*}
     \Mklz{ I_C\otimes \FK{D}\,+\,\FK{C}\otimes I_D }{ I_C\in \FC,\;I_D\in \FD}.
  \end{eqnarray*}
This gives a F-variety $(C\times D, \FK{C\times D}, \Fi{C\times D})$ which is, together with the projections $pr_C$, $pr_D$, the 
product of $(C,\FK{C},\Fi{C})\,$, $(D,\FK{D},\Fi{D})$ in the category of F-varieties.
Let $c\in C$ and $d\in D$. The map
\begin{eqnarray*} 
   T_{(c,d)}\, pr_C\times T_{(c,d)}\, pr_D:\;T_{(c,d)}(C\times D) \;\to\; T_c C \times T_d D
\end{eqnarray*} 
is bijective.\vspace*{1ex} 
%

{\bf Principal open sets:} Let $(A,\FK{A},\Fi{A})$ be a F-variety, and let $g\in\FK{A}\setminus \{0\}$.  
Equip the coordinate ring $\FK{D(g)}\cong \FK{A}_g$ of the principal open set
$D(g)$ with the filter of ideals
\begin{eqnarray*}
    \Mklz{ \,I^e\,}{\,I\in\FA}\;,&\mb{ where }& I^e:\;=\;\Mklz{ \frac{f}{g^n} }{ f\in I\;,\;n\in \Nn }. 
\end{eqnarray*} 
This gives a F-variety $(D(g),\FK{D(g)},\Fi{D(g)})$.
The inclusion map $j:D(g)\to A$ is a morphism, and for $a\in D(g)$ the tangent map  
$ T_a j:\,T_a (D(g))\to T_a A $ is bijective.\vspace*{1ex} 
%
%

{\bf The category of weak F-algebraic monoids and groups:} Let $M$ be a monoid, and let $(M,\FK{M},\Fi{M})$ be a variety, such that for
all $m\in M$ the left and right multiplications $r_m$, 
$l_m$ by $m$ are morphisms of F-varieties. We get an injective linear map 
\begin{eqnarray*}
  \Psi_l:\, T_1 M &\to & Der(\FK{M})
\end{eqnarray*}
by $(\Psi_l(\delta)f)(m):= \delta(l_m^* f)$, $f\in\FK{M}$, $m\in M$. Its image
consists of left invariant derivations of $\FK{M}$. 
Similarly, we get an injective linear map 
\begin{eqnarray*}
  \Psi_r:\, T_1 M &\to & Der(\FK{M})
\end{eqnarray*}
by $(\Psi_r(\delta)f)(m):= \delta(r_m^* f)$, $f\in\FK{M}$, $m\in M$. Its image
consists of right invariant derivations of $\FK{M}$. 
The inverse maps on the images of $\Psi_l$ and $\Psi_r$ are obtained by
concatenation with the evaluation morphism of $\FK{M}$ in $1\in M$. 
We call $(M,\FK{M},\Fi{M})$ a {\it weak F-algebraic monoid}, if the following
things hold:
\begin{itemize}
\item The images of $\Psi_{l}$ and $\Psi_{r}$ are 
Lie subalgebras of $Der(\FK{M})$, and the Lie algebra structures on $T_1M$, obtained by pulling back, are opposite.
The tangent space $T_1 M$ equipped with the structure of the Lie algebra, which is obtained by pulling back by $\Psi_{l}$, 
is called the {\it Lie algebra} $Lie(M)$ of $M$.
\item For every $g\in M^\times$ and $\delta\in Lie(M)$ also
  $g.\delta:=\delta\circ l_g^*\circ r_{g^{-1}}^*\in Lie(M)$. This
  defines an action of $M^\times$ on $Lie(M)$, which we call the {\it adjoint action}.
\end{itemize}
A {\it morphism} of weak F-algebraic monoids $(M,\FK{M},\Fi{M})$, $(N,\FK{N},\Fi{N})$ consists of a homomorphism of monoids 
$\phi:M\to N$, which is also a morphism of F-varieties. Then the map 
  \begin{eqnarray*}
     Lie(\phi)\,:=\,T_1\phi: \;Lie(M) &\to & Lie(N)
  \end{eqnarray*}
is a homomorphism of Lie algebras.\vspace*{1ex}

A weak F-algebraic monoid $(M,\FK{M},\Fi{M})$ {\it acts algebraically} on the F-variety
$(A,\FK{A},\Fi{A})$ if the monoid $M$ acts on the set $A$ such that for all
$m\in M$ the left applications $l_m:A\to A$ defined by $l_m(a):=ma$, $a\in A$,
and for all $a\in A$ the evaluations $r_a:M\to A$ defined by $r_a(m):=ma$,
$m\in M$, are morphisms of F-varieties.\vspace*{1ex}

Let $(M,\FK{M},\Fi{M})$ be a weak F-algebraic monoid. Every closed
submonoid $N$ of $M$, equipped with its subvariety structure $(N,\FK{N},\Fi{N})$, is again a weak
F-algebraic monoid. If $(M,\FK{M},\Fi{M})$ acts algebraically on
$(A,\FK{A},\Fi{A})$, then also $(N,\FK{N},\Fi{N})$ does.\vspace*{1ex}

Let $(M,\FK{M},\Fi{M})$ be a weak F-algebraic monoid with an involution of monoids $*:M\to M$, which is also an isomorphism of 
F-varieties. This involution induces an involution $*$ of the algebra $\FK{M}$
by $f^*(m):=f(m^*)$, $f\in\FK{M}$, $m\in M$. It is easy to see that we 
also get an involution $*$ of the Lie algebra $Lie(M)$ by $\delta^*:=\delta\circ *\,$, $\delta\in Lie(M)$.

A {\it weak F-algebraic group}  consists of a group $G$ equipped with the
structure of a weak F-algebraic monoid $(G,\FK{G},\Fi{G})$, such that the
inverse map $\mb{}^{-1}:G\to G$ is an isomorphism of F-varieties, and the corresponding
involution on the Lie algebra $Lie(M)$ is the multiplication by $-1\in\F$.\vspace*{1ex}
%
%
%

{\bf Some remarks to F-varieties:} The full subcategory of the category of F-varieties given by the F-varieties with discrete filter of ideals 
coincides with  the full subcategory of the category of sets with coordinate rings given by the objects $(A,\FK{A})$ such that, by our identification, 
$A=\Spm\FK{A}$. Here $T_a A=Der_a(A)$ for all $a\in A$.\vspace*{1ex}

The algebraic geometry of ind-varieties and ind-groups has been developed in
\cite{Sa}. See also Chapter IV of the book \cite{Ku}.
The category of affine ind-varieties identifies with a full subcategory of the category of F-varieties. The tangent spaces of an affine 
ind-variety and its corresponding F-variety coincide. (But taking
ind-subvarieties and F-subvarieties is slightly different.)

We do not explain all details here, only how to make and affine ind-variety into an F-variety: 
An affine ind-variety structure on a non-empty set $X$ is given by a family $\mathcal A$ of subsets of $X$, each equipped with the structure 
of an affine algebraic variety, such that the following properties hold:
The sets of $\mathcal A$ cover $X$. The system of sets $\mathcal A$ is directed upwards, i.e., for all $A_1,\,A_2\in{\mathcal A}$ 
there exits a set $A\in {\mathcal A}$ such that $A_1,\,A_2\subseteq A$. If $A_1,\,A_2\in {\mathcal A}$ and $A_1\subseteq A_2$,
then $A_1$ is closed in $A_2$, and its variety structure is the closed subvariety structure induced by $A_2$.
There may be several defining systems $\mathcal A$ for the same ind-variety
structure, but all are included in a biggest defining system, the
complete defining system.
The set $X$ is the inductive limit of the subsets given by $\mathcal A$. The coordinate ring $\FK{X}$ of the affine ind-variety $X$ with its standard 
structure of a topological algebra is given by the inverse limit of the discrete topological algebras $\FK{A}$, $A\in {\mathcal A}$.
The open ideals of $\FK{X}$ are a topological zero filterbase, in particular they determine the topological algebra structure on $\FK{X}$.

The vanishing ideals $I(A)$, $A\in{\mathcal A}$,  give a filterbase ${\mathcal B}_X$ of ideals of $\FK{X}$. 
The corresponding filter of ideals $\Fi{X}$ coincides with the set of open
ideals of $\FK{X}$. In this way the affine ind-variety $X$ is also a 
F-variety. 
We can recover the ind-variety structure from its F-variety structure. 
The zero sets of the ideals of ${\mathcal B}_X$ coincide
  with the sets of ${\mathcal A}$, the coordinate ring of such a set is
  obtained by restricting the coordinate ring of $X$. 
Similarly, the zero sets of the ideals of $\Fi{X}$ determine the complete defining system of $X$.

\subsection{Varieties and weak algebraic monoids and groups\label{var}}
In elementary algebraic geometry affine algebraic varieties can be constructed
as follows: Equip the finite-dimensional affine spaces with their canonical
coordinate rings. Equip the closed subsets with their coordinate rings obtained
by restriction. An affine algebraic variety is a set with coordinate ring
which is isomorphic to some closed subset of some finite
dimensional affine space.
The varieties which we need in this paper are constructed similarly.\vspace*{1ex}

Let $V$ be a $\F$-linear space. For a subspace $V^{du}$ of the full dual $V^*$,
which separates the points of $V$, we get a F-variety $(V,\FK{V},\Fi{V})$ in
the following way:
\begin{itemize}
\item The coordinate ring $\FK{V}$ is the algebra of functions generated by
  $V^{du}$. It is isomorphic to the symmetric algebra in $V^{du}$ in the obvious way. 
\item The filter of ideals $\Fi{V}$ is determined by the filterbase  
  \begin{eqnarray*}
    \Mklz{ I_V(U)}{ U \mb{ a finite dimensional subspace of } V },
  \end{eqnarray*}
where $I_V(U)$ denotes the ideal of functions of $\FK{V}$ vanishing on $U$.
\end{itemize}

The tangent space $T_a V$ at $a\in V$ can be canonically identified with $V$ by means of the linear
bijective map $V\to T_a V\,$, which assigns to $v\in V$ the derivation $\delta_v\in T_a V$ defined by
  \begin{eqnarray*}
    f(a+tv)=f(a)+t\delta_v (f) + \mb{ higher terms in }t\quad\mb{ where }\quad  t\in \F,\;f\in\FK{V}.
  \end{eqnarray*}
Equivalently, $\delta_v \in T_a V$ is the derivation of $\FK{V}$ in $a$
determined by $\delta_v(\phi)=\phi(v)$, $\phi\in V^{du}$.

\begin{Remark}
If $V$ is of infinite dimension there are many possibilities to choose a subspace $V^{du}$ of $V^*$, which separates 
the points of $V$. Therefore there are many possible coordinate rings of $V$.
The situation is different if we restrict to a finite dimensional subspace $U$ of $V$. Because $V^{du}\res{U}$ separates the points of 
$U$, there is only the possibility $V^{du}\res{U}\,=U^*$, and $\FK{U}$ is the classical coordinate ring of $U$.
The induced filter of ideals of $\FK{U}$ is discrete because of $I(U)\res{U}\,=\{0\}\in\Fi{U}$.

Independently of the chosen point separating subspace $V^{du}\subseteq V^*$, the
tangent spaces $T_a V$, $a\in V$, identify with $V$.
\end{Remark}
\begin{Remark} \label{UVR}Let $V$ be a $\F$-linear space, let $V^{du}$ be a subspace of
  $V^*$, which separates the points of $V$, and let $(V,\FK{V},\Fi{V})$ be the
  corresponding  F-variety. Let $W$ be a closed $\F$-linear subspace of $V$.
  The restriction $W^{du}$ of $V^{du}$ onto $W$ is a subspace of $W^*$, which
  is point separating on $W$. The F-variety structure given
  by $W^{du}$ coincides with the F-subvariety structure on $W$ induced by $(V,\FK{V},\Fi{V})$. 
\end{Remark}
\begin{Definition}\mb{}\label{V1}
A variety is a F-variety isomorphic to a subvariety of $(V,\FK{V},\Fi{V})$ for some
$\F$-linear space $V$, and point separating subspace $V^{du}\subseteq V^*$. Morphisms of varieties are morphism of F-varieties.
\end{Definition}
\begin{Remarks} 
The category of varieties includes the category of the classical affine
algebraic varieties as a full subcategory, if we equip the coordinate rings
of the affine algebraic varieties with the discrete filter of ideals.
\end{Remarks}

It is not difficult to show:
\begin{Proposition}\label{V2} Closed subsets, products, and principal open sets of varieties are again varieties. 
\end{Proposition}
\begin{Definition}\label{V3} A weak algebraic monoid resp. group is a weak F-algebraic monoid resp. group, whose underlying F-variety is a variety.
\end{Definition}

There is the following easy possibility to obtain weak algebraic monoids: 
Let $A$ be an associative $\F$-algebra with unit. Let $A^{du}$ be a subspace
of $A^*$, which separates the points of $A$, such that for all $a\in A$ there exist the
adjoints $l_a^{du}:A^{du}\to A^{du}$ and $r_a^{du}:A^{du}\to A^{du}$ of the
left and right multiplications $l_a$ and $r_a$ by $a$. 
Equip $A$ with its canonical structure of a variety $(A,\FK{A},\Fi{A})$
induced by $A^{du}$. Then $A$ is a weak algebraic monoid. By the canonical identification
of $T_1 A$ with $A$, the Lie algebra
$Lie(A)$ identifies with the Lie algebra $A$ associated to the associative
algebra $A$.
There are two important examples:\vspace*{1ex}

(1) Let $V$ be a $\F$-linear space. Let $V^{du}$ be a subspace of $V^*$ which separates the points of $V$.   
Recall that $End_{V^{du}}(V)$ denotes the subalgebra of $End(V)$ 
consisting of the endomorphisms $\al:\, V\to V$, for which the dual map $\al^{du}:\,
V^{du}\to V^{du}$ exists. The functions
\begin{eqnarray*}
\begin{array}{cccc}
  \ti{f}_{\phi v}: &End_{V^{du}}(V)&\to & \F\\
                   &  \al &\mapsto & \phi(\al v )
\end{array} \quad\mb{ where }\quad v\in V,\;\phi\in V^{du}
\end{eqnarray*} 
span a subspace $(End_{V^{du}}(V))^{du}\cong V^{du}\otimes V$ of
$(End_{V^{du}}(V))^*$ which has the required properties. $End_{V^{du}}(V)$
equipped with the corresponding variety structure 
is a weak algebraic monoid. 
Its Lie algebra $Lie(End_{V^{du}}(V))$ identifies canonically with the Lie algebra $End_{V^{du}}(V)$.

Equip $V$ with its variety structure induced by $V^{du}$. Then $End_{V^{du}}(V)$ acts on $V$ algebraically.\vspace*{1ex}

(2) Let $\mathcal C$ be a category of representations of a Lie algebra $\g$ as
in Section \ref{TK}, and ${\mathcal C}^{du}$ a category of duals.
Recall from Section \ref{TK} the $\F$-algebra $Nat$. For an object $V$ of
$\cal C$, $\phi\in V^{du}$, and $v\in V$ define a linear function $\ti{f}_{\phi v}: Nat\to\F$ by
\begin{eqnarray*}
   \ti{f}_{\phi v}(m) :=  \phi(m_V v)\quad,\quad m\in Nat .
\end{eqnarray*}
Then by equation (\ref{msum}) the set
\begin{eqnarray*}
      Nat^{du}:=\Mklz{\ti{f}_{\phi v}\,}{\,\phi\in V^{du}\,,\, v\in V\,,\,\,V \;an\; object\;of\; {\mathcal C} }.
\end{eqnarray*}
is a $\F$-linear subspace of $Nat^*$, which has the required properties. $Nat$
equipped with the
corresponding variety structure $(Nat,\FK{Nat},\Fi{Nat})$ is a weak algebraic
monoid. Its Lie algebra $Lie(Nat)$ identifies canonically with the Lie algebra $Nat$.

Let $V$ be an object of $\mathcal C$. Equip $V$ with its variety structure induced by $V^{du}$. Then $Nat$ acts on $V$ algebraically.\vspace*{1ex}
\begin{Theorem}\label{V4} Let $M$ be the Tannaka monoid with coordinate ring of matrix coefficients $\FK{M}$
  associated to the categories $\mathcal C$, $\mathcal C^{du}$ by the Tannaka
  reconstruction of Section \ref{TK}. Let $Lie(M)$ be its Lie algebra as defined in Section \ref{TK}. Then:
\begin{itemize}
\item[(a)] The Tannaka monoid $M$ is a closed submonoid of $Nat$.
\item[(b)] The coordinate ring $\FK{M}$ coincides with the restriction of the coordinate ring $\FK{Nat}$ onto $M$.
\item[(c)] Assigning to an element of $Lie(M)$ its corresponding derivation of $\FK{M}$ in $1\in M$, the Lie algebra $Lie(M)$ identifies with the Lie algebra of the weak
  algebraic monoid $(M,\FK{M},\Fi{M})$, where $\Fi{M}$ is the canonical restriction of the filter of ideals $\Fi{Nat}$.
\end{itemize}
\end{Theorem}
\begin{Remark} The Lie algebra of the weak algebraic monoid $(Nat,\FK{Nat},\Fi{Nat})$ identifies canonically with the Lie algebra $Nat$. The Lie algebra of the
weak algebraic monoid $(N,\FK{N},\Fi{N})$, given by a submonoid $N\subseteq Nat$, identifies with a Lie subalgebra of $Nat$.
In Section \ref{TK} our Lie algebras for submonoids of $M$ have been defined in this way, as submonoids of $Nat$.
\end{Remark}

\Proof
To (a): Recall that the monoid $M$ has been defined as the submonoid of $Nat$,
given by the elements $m\in Nat$, which satisfy 
\begin{itemize}
\item[(1)] $\quad m_{V\otimes W}\,=\, m_V\otimes m_W\;$ for all objects $V$, $W$,
\item[(2)] $\quad m_{V_0}\, =\, id_{V_0}\,$ for every trivial one-dimensional $\g$-module.
\end{itemize}
We have to show $\overline{M}\subseteq M$.
Due to (1), for all $m\in M$, and for all objects $V$, $W$ of $\mathcal C$, for all $v\in V$, 
$w\in W$, $\phi\in V^{du}$, and $\psi\in W^{du}$ we get
\begin{eqnarray}\label{mm1}
   \ti{f}_{\phi\otimes\psi\,v\otimes w}(m)- \ti{f}_{\phi\,v}(m)\ti{f}_{\psi\,w}(m)=0.
\end{eqnarray}
Due (2), for all $m\in M$, and all $v_0\in V_0$ and $\phi_0\in V_0^{du}$ we get
\begin{eqnarray}\label{mm2}
     \ti{f}_{\phi_0\,v_0}(m)-\phi_0(v_0)=0.
\end{eqnarray}
Now let $\overline{m}\in\overline{M}$. Since the equations (\ref{mm1}),
(\ref{mm2}) are valid for all $m\in M$ they are also 
valid for $\overline{m}$. Since $V^{du}$, $W^{du}$, and $V_0^{du}$ are point
separating, this implies the properties (1), (2) for $\overline{m}$. 

Part (b) of the theorem is obvious by the definition of the coordinate ring
$\FK{M}$. To (c): Recall that the Lie algebra $Lie(M)$ can be defined as the
Lie subalgebra of $Nat$ given by the elements $x\in Nat$ for which
there exists a derivation $\delta_x:\FK{M}\to \F$ in $1\in M$ such that 
\begin{eqnarray}\label{dprop}
   \delta_x(f_{\phi v})\;=\;\phi(x_V v)\; \mb{ for all } \; \phi\in V^{du}\,,\;v\in V\,, \;V \mb{ an object of } {\mathcal C}.
\end{eqnarray} 
Denote by $\iota:M\to Nat$ the
inclusion map. For this proof we denote the Lie
algebra of the weak algebraic monoid $(M,\FK{M},\Fi{M})$ by $T_1 M$, and the
Lie algebra of the weak algebraic monoid $(Nat,\FK{Nat},\Fi{Nat})$ by $T_1
Nat$. The embedding $\iota:M\to Nat$ gives an embedding of Lie algebras
$T_1\iota:T_1 M\to T_1 Nat$ with image 
\begin{eqnarray*}
    T_1\iota (T_1 M)= \Mklz{\delta\in T_1 Nat}{\delta (I(M))=0}\subseteq T_1 Nat,
\end{eqnarray*}
where $I(M)$ denotes the vanishing ideal of $M$ in $\FK{Nat}$. By the
canonical identification of $Nat$ with $T_1 Nat$, mapping $x\in Nat$ to the
derivation $\ti{\delta}_x\in T_1 Nat$, the Lie algebra $T_1\iota (T_1 M)$
identifies with the Lie subalgebra
\begin{eqnarray*}
    \Mklz{ x \in Nat}{\ti{\delta}_x (I(M))=0}
\end{eqnarray*}
of $Nat$. If $x\in Nat$, then by the definition of $\ti{\delta}_x$ we have
$\ti{\delta}_x (I(M))=0$ if and only if
$\ti{\delta}_x$ factors to a derivation $\delta_x:\FK{M}\to \F$ with property (\ref{dprop}).
\qed

Let $V$ be a $\F$-linear space. Let $(V,\FK{V},\Fi{V})$ be the variety induced by a subspace $V^{du}\subseteq V^*$ which separates the 
points of $V$. Identify the tangent spaces $T_v V$, $v\in V$, canonically with $V$. Let $(M,\FK{M},\Fi{M})$ be a weak F-algebraic monoid 
which acts algebraically and linearly on $(V,\FK{V},\Fi{V})$. Recall that for $v\in V$ the map $r_v:M\to V$ 
is defined by $r_v(m):=mv$, $m\in M$. It is not difficult to check that we get an action of the Lie algebra $Lie(M)$ of $M$, 
which we call the differentiated action, by
\begin{eqnarray*}
   \delta.v:=(T_1 r_v)(\delta)\in T_v V=V \quad\mb{ where }\quad \delta\in Lie(M),\;v\in V.
\end{eqnarray*}  

\begin{Corollary} Let $(M,\FK{M},\Fi{M})$ be the weak algebraic monoid associated to the categories $\mathcal C$, $\mathcal C^{du}$ 
by the Tannaka reconstruction of Section \ref{TK}. Identify the Lie algebra $Lie(M)$ as defined in Section \ref{TK} with the Lie algebra 
of the weak algebraic monoid $(M,\FK{M},\Fi{M})$.
Let $V$ be an object of $\mathcal C$ equipped with the variety structure $(V,\FK{V},\Fi{V})$ induced by $V^{du}$. Then:
\begin{itemize}
\item[(a)] The weak algebraic monoid $(M,\FK{M},\Fi{M})$ acts algebraically
  and linearly on $(V,\FK{V},\Fi{V})$.
\item[(b)] The action of $Lie(M)$ on $V$ of Section \ref{TK}
  coincides with the differentiated action.\vspace*{1ex}
\end{itemize}
\end{Corollary}

Let $\mathcal C$ be a category of representations of a Lie algebra $\g$ as in Section \ref{TK}, and ${\mathcal C}^{du}$ a category of duals.
Sometimes it is useful to identify the variety $Nat$ with varieties $Nat({\mathcal G})$ associated to certain full subcategories 
$\mathcal G$ and the corresponding categories of duals ${\mathcal G}^{du}$.

We only treat a particular case for semisimple categories $\mathcal C$, which will be important in \cite{M4}.  
Recall the notation introduced to state Theorem \ref{PW}. 
Equip the algebra 
\begin{eqnarray*}
    \prod_{W\in Irr } End_{W^{du}}(W)
\end{eqnarray*}
with the structure of a variety induced by 
\begin{eqnarray*}
   \left(\prod_{W\in Irr}
     End_{W^{du}}(W)\right)^{du}:=\bigoplus_{W\in {\mathcal G}}
     \left( End_{W^{du}}(W)\right)^{du} \subseteq 
                                       \left(\prod_{W\in {\mathcal G}} End_{W^{du}}(W)\right)^*.
\end{eqnarray*}
If $V\in Irr $ denote by $\left(End_{V^{du}}(V)\right)_{S(V)}$ the $S(V)$-linear endomorphisms of $End_{V^{du}}(V)$. 
It is easy to check that 
\begin{eqnarray*}
   Nat(Irr):=\prod_{W\in Irr } \left(End_{W^{du}}(W)\right)_{S(W)}
\end{eqnarray*} 
is a closed $\F$-subalgebra of $\prod_{W\in Irr } End_{W^{du}}(W)$, which we equip with its subvariety structure.  

\begin{Theorem}\label{IrrM} Let $\mathcal C$ be semisimple. Suppose that for all non-empty families $(V_j )_{j\in J}$ of elements of $Irr$, 
such that $\oplus_{j\in J}V_j$ is an object of $\mathcal C$, and for all families $(m_j)_{j\in J}$ of endomorphisms 
$m_j\in \left(End_{(V_j)^{du}}(V_j)\right)_{S(V_j)}$ the endomorphism $\oplus_{j\in J} m_j$ is already in 
$End_{(\oplus_j V_j)^{du}}(\oplus_j V_j)$. Then
\begin{eqnarray*} 
  \Xi_{Irr}:\;Nat \qquad\quad &\to&  \;Nat(Irr)\\
   m=(m_V)_{V\;an\; obj.\;of\;{\mathcal C}} &\mapsto &(m_W)_{W\in Irr}
\end{eqnarray*}
is an isomorphism of $\F$-algebras, which is also an isomorphism of varieties.
\end{Theorem}

\begin{Remark}The conditions of the theorem on the endomorphisms are satisfied
  if $V^{du}=V^*$ for all objects of $\mathcal C$. They are also satisfied if
  for all non-empty families $(V_j )_{j\in J}$ of elements of $Irr$, such that
  $\oplus_{j\in J}V_j$ is an object of $\mathcal C$, we have 
\begin{eqnarray*}
   (\,\bigoplus_{j\in J}V_j\,)^{du}=\bigoplus_{j\in J} \,(V_j)^{du}.
\end{eqnarray*}
\end{Remark}

\Proof  It is easy to check that $\Xi:=\Xi_{Irr}$ is a morphism of $\F$-algebras, whose comorphism $\Xi^*$ exists and maps 
$Nat(Irr)^{du}=\left(\prod_{W\in Irr} End_{W^{du}}(W)\right)^{du}\res{Nat(Irr)}$ onto $Nat^{du}$. 
The Theorem follows if we show that $\Xi$ is bijective. We construct a map $\Omega$ which is inverse to $\Xi$. 

(a) Let $n=(n_V)_{V\in Irr}\in Nat(Irr)$. Let $\bigoplus_{p\in P} V_p$ and $\bigoplus_{q\in Q}W_q$ be two sums of
$\g$-modules of $Irr$. For $\ti{p}\in P$ let $in_{\ti{p}}$ be the canonical injection of $V_{\ti{p}}$ into $\bigoplus_{p\in P} V_p$.
For $\ti{q}\in Q$ let $pr_{\ti{q}}$ be the canonical projection of $\bigoplus_{q\in Q} W_q$ onto $W_{\ti{q}}$. 
Set $n_P:=\bigoplus_{p\in P}n_{V_p}$ and $n_Q:=\bigoplus_{q\in Q}n_{W_q}$. Let
$\gamma: \,\bigoplus_{p\in P} V_p\to \bigoplus_{q\in Q}W_q$ be a homomorphism
of $\g$-modules. Then
\begin{eqnarray*}
  \gamma\circ n_P=n_Q\circ \gamma
\end{eqnarray*} 
if and only if 
\begin{eqnarray}\label{inpr}
  (pr_q\circ\gamma\circ in_p)\circ n_{V_p}= n_{W_q}\circ (pr_q\circ\gamma\circ in_p ) 
\end{eqnarray}
for all $p\in P$, $q\in Q$. The map $pr_q\circ\gamma\circ in_p :\,V_p\to W_q$ is a homomorphism between
the irreducible $\g$-modules $V_p,\,W_q\in Irr$. If $V_p\neq W_q$ then
(\ref{inpr}) is valid since $pr_q\circ\gamma\circ in_p$ is the zero map. If $V_p=W_q$ then
(\ref{inpr}) is valid since $n_V\in (End_{V^{du}}(V))_{S(V)}$, $V\in Irr$.

(b) Let $n=(n_V)_{V\in Irr}\in Nat(Irr)$. For every object $V$ of $\mathcal C$ fix an isomorphism 
$\al_V:\bigoplus_{p\in P(V)} V_p \to V $ of $\g$-modules. Define
\begin{eqnarray*}
  n_V :=\al_V \circ n_{P(V)}\circ \al_V^{-1}.
\end{eqnarray*}
We have $n_V\in End_{V^{du}}(V)$ since the adjoint maps of $\al_V$, $n_{P(V)}$,
and $\al_V^{-1}$ exist. This definition is independent of the chosen
isomorphism. Let $\beta: \bigoplus_{p\in P(V)} V_p\to V $ be another isomorphism. Then  
$\al_V \circ n_{P(V)}\circ \al_V^{-1}=\beta\circ n_{P(V)}\circ \beta^{-1}$
if and only if 
$(\beta^{-1}\circ\al_V )\circ n_{P(V)}= n_{P(V)}\circ (\beta^{-1}\circ\al_V)$,
which is valid by (a).

Let $V$, $W$ be objects of $\mathcal C$. Let $\beta\in Hom_{\bf g}(V,W)$. Then by definition 
$\beta\circ n_V=n_W\circ \beta$ 
if and only if 
$\beta\circ\al_V\circ n_{P(V)}\circ\al_V^{-1}=\al_W \circ n_{P(W)}\circ\al_W^{-1}\circ\beta$ 
if and only if 
$(\al_W^{-1}\circ\beta\circ\al_V)\circ n_{P(V)}=n_{P(W)}\circ(\al_W^{-1}\circ\beta\circ\al_V)$, 
which is valid by (a).

Therefore $\Omega(n):= (n_V)_{V \,an\, object\, of\, {\mathcal C} } \in Nat$. The map $\Omega : Nat(Irr)\to Nat$ satisfies 
$\Xi\circ\Omega=id_{Nat(Irr)}$ and $\Omega\circ \Xi=id_{Nat}$. 
\qed

{\bf Some remarks to varieties:} 
A more general class of varieties can be obtained from the F-varieties $(V,\FK{V},\Fi{V})$ of $\F$-linear spaces $V$, where:
\begin{itemize}
\item[(a)] The coordinate ring $\FK{V}$ is an algebra of functions, which restricts
on every finite-dimensional linear subspace $U$ of $V$ to the classical coordinate ring of $U$. 
\item[(b)] $\FK{V}\cap V^*$ restricts on every finite-dimensional linear
  subspace $U$ of $V$ to $U^*$. Equivalently $\FK{V}\cap V^*$ is point
  separating on $V$. 
\item[(c)] The filter of ideals $\Fi{V}$ is obtained as above by the filterbase
  of vanishing ideals $I_V(U)$, $U$ a finite-dimensional linear subspace of $V$.
\end{itemize}
(Maybe (a) and (b) are related.) Note that the coordinate ring $\FK{V}$ can be looked at as a completion of the coordinate ring
associated to $V^{du}:=\FK{V}\cap V^*$.
Now the F-varieties associated to the ind-varieties of $\F$-linear spaces are
included as varieties. For such a variety $V^{du}=V^*$.\vspace*{1ex}

Parallel, it is possible to generalize the Tannaka-Krein reconstruction to categories $\mathcal C$ of $\g$-modules as above, 
where in addition the $\g$-modules contained in $\mathcal C$ are equipped with coordinate rings of this sort, satisfying certain
compatibility conditions. Note also that $V^{du}=\FK{V}\cap V^*$, $V$ an
object of $\mathcal C$, gives a category of duals ${\mathcal C}^{du}$. \vspace*{1ex} 

In the situation which we consider in this paper, choosing a category of duals, the coordinate rings of the monoids are directly related to 
representation theory. In addition, the coordinate rings of the $\F$-linear spaces are, similarly as in the classical case, symmetric algebras. 
In this more general situation the monoids are submonoids of the monoids
obtained by the corresponding category of duals. The coordinate rings of the
monoids can be looked at as completions of the coordinate rings obtained by
the corresponding category of duals. 
%
%
%
%

%
%
%

%
%
\end{document}